\documentclass[a4paper]{amsart}
\usepackage[dvips]{graphicx}
\theoremstyle{plain}
  \newtheorem{thm}{Theorem}[section]
  \newtheorem{lem}[thm]{Lemma}
  
  \newtheorem{prop}[thm]{Proposition}
  \newtheorem{conj}[thm]{Conjecture}

\theoremstyle{definition}

\theoremstyle{remark}
  \newtheorem{rem}[thm]{Remark}
  \newtheorem*{ack}{Acknowledgments}

\newcommand{\C}{\mathbb{C}}
\newcommand{\Vol}{\operatorname{Vol}}
\numberwithin{equation}{section}
\numberwithin{figure}{section}
\date{\today}
\allowdisplaybreaks
\begin{document}
\title[Kashaev's invariant and the volume of a hyperbolic knot]
{Kashaev's invariant and the volume of a hyperbolic knot after Y.~Yokota}
\author{Hitoshi Murakami}
\address{
  Department of Mathematics,
  Tokyo Institute of Technology,
  Oh-okayama, Meguro, Tokyo 152-8551, Japan
}
\email{starshea@tky3.3web.ne.jp}
\thanks{
This research is supported in part by Grand-in-Aid for Scientific Research,
Japan Society for the Promotion of Science.
}
\begin{abstract}
I follow Y.~Yokota to explain how to obtain a tetrahedron decomposition
of the complement of a hyperbolic knot and compare it with the asymptotic
behavior of Kashaev's link invariant using the figure-eight knot as an example.
\end{abstract}
\maketitle
\section{Introduction}\label{sec:introduction}
In \cite{Kashaev:MODPLA95} R.~Kashaev introduced a link invariant by using
the quantum dilogarithm and after some concrete calculations for the knots $4_1$,
$5_2$, and $6_1$, he conjectured
in \cite{Kashaev:LETMP97} the asymptotic behavior of his invariant determines
the hyperbolic volume for hyperbolic links.
More precisely his conjecture is as follows.
For an integer $N\ge2$ let $\langle{L}\rangle_N$ be the Kashaev invariant of
a link $L$.
\begin{conj}[Kashaev]
  If $L$ is hyperbolic, i.e. $S^3\setminus{L}$ has a complete hyperbolic
  structure, then
  \begin{equation*}
    2\pi\lim_{N\to\infty}\frac{\log\left|\langle{L}\rangle_N\right|}{N}=\Vol(L),
  \end{equation*}
  where $\Vol(L)$ is the volume of $S^3\setminus{L}$.
\end{conj}
\par
On the other hand in \cite{Murakami/Murakami:volume}, J.~Murakami and I proved
that Kashaev's invariant coincides with a colored Jones polynomial evaluated
at a root of unity.
By using this Y.~Yokota, following D.~Thurston \cite{D.Thurston:Grenoble},
considered a tetrahedron decomposition corresponding to the $R$-matrix used
to define colored Jones polynomials and showed that at least for a concrete
example, the knot $6_2$, this decomposition gives a nice and constructive proof
for Kashaev's conjecture \cite{Yokota:Murasugi70}
(he also mentioned that this seems work well for other knots; especially for
alternating knots).
Several computations (without tetrahedron decompositions) were made in
\cite{Murakami/Murakami/Okamoto/Takata/Yokota:CS} and the conjecture was also
confirmed for the knots $6_3$, $7_2$, and $8_9$, and for the Whitehead link.
(In \cite{Murakami/Murakami/Okamoto/Takata/Yokota:CS} a relation between the
asymptotic behavior of the Kashaev invariant and the Chern--Simons invariant
is also mentioned.)
\par
After that Yokota used Kashaev's original $R$-matrix and found another
tetrahedron decomposition of a hyperbolic knot complement fit for
the asymptotic behavior of Kashaev's invariant \cite{Yokota:volume}
(see also \cite{Yokota:Topology_Symposium2000}).
This gives simpler decompositions and seems to work much better including
non-alternating knots.
\par
This is an expository note to describe Yokota's method explaining beautiful
harmony of Kashaev's invariants and hyperbolic structures, putting special
emphasis on the figure-eight knot.
\begin{ack}
Part of this work was prepared for a series of lectures at Chiba University
in July, 2000.
I would like to thank K.~Kuga for giving me the opportunity to give lectures
there.
Thanks are also due to J.~Murakami, M.~Okamoto, T.~Takata and Y.~Yokota for
helpful discussions.
\end{ack}
\section{Algebra}\label{sec:algebra}
In this section I describe Kashaev's link invariant \cite{Kashaev:MODPLA95}
following Yokota.
We fix an integer $N\ge2$ and put $q=\exp\left(2\pi\sqrt{-1}/N\right)$.
Put $(x)_n=\prod_{k=1}^{n}(1-x^k)$ for $n\ge 0$.
For integers $0{\le}k,l,m,n{\le}N-1$ we put
\begin{equation*}
  \theta\!\left(\begin{smallmatrix} k & n \\ l & m \end{smallmatrix}\right)
  =
  \begin{cases}
    1\;&\text{if ${l}\ne{m}$ and $q^{k},q^{l},q^{m},q^{n}$ are on the unit
    circle counterclockwise}
    \\
    &\hfill\text{(other pairs may be the same)},
    \\
    0\;&\text{otherwise}.
  \end{cases}
\end{equation*}
Note that
$\theta\!\left(\begin{smallmatrix} k & n \\ l & m \end{smallmatrix}\right)=1$
if and only if
\begin{equation*}
  {k}\le{l}<{m}\le{n},\quad
  {n}\le{k}\le{l}<{m},\quad
  {m}\le{n}\le{k}\le{l}\;(m<l),\;\text{or}\quad
  {l}<{m}\le{n}\le{k}.
\end{equation*}
For an integer $x$, we denote by $[x]\in\{0,1,2,\dots,N-1\}$ the residue
modulo $N$.
\par
Now Kashaev's $R$-matrix $R=\left(R_{l,m}^{k,n}\right)$ is an
$N^2\times N^2$-matrix given by
\begin{align*}
  R_{l,m}^{k,n}
  =
  \frac
  {Nq^{1-(l-n+1)(m-k)}
  \theta\!\left(\begin{smallmatrix} k & n \\ l & m \end{smallmatrix}\right)}
  {(q)_{[m-l-1]}(q^{-1})_{[n-m]}(q)_{[k-n]}(q^{-1})_{[l-k]}}
  \\ \intertext{and}
  \left(R^{-1}\right)_{l,m}^{k,n}
  =
  \frac
  {Nq^{-1+(m-k-1)(l-n)}
  \theta\!\left(\begin{smallmatrix} k & n \\ l & m \end{smallmatrix}\right)}
  {(q^{-1})_{[m-l-1]}(q)_{[n-m]}(q^{-1})_{[k-n]}(q)_{[l-k]}}.
\end{align*}
Let $\mu$ be an $N\times N$-matrix with $\mu_l^k=-\delta_{k,l+1}q^{1/2}$,
where $\delta_{i,j}$ is Kronecker's delta.
Then $(R,\mu,-q^{1/2},1)$ is an enhanced Yang--Baxter operator
\cite{Turaev:INVEM88}, i.e. the following equalities hold.
(See \cite{Murakami/Murakami:volume} for a proof.)
\begin{prop}
  \begin{gather*}
    \left({R}\otimes{id}_{\C^N}\right)
    \left({id}_{\C^N}\otimes{R}\right)
    \left({R}\otimes{id}_{\C^N}\right)
    =
    \left({id}_{\C^N}\otimes{R}\right)
    \left({R}\otimes{id}_{\C^N}\right)
    \left({id}_{\C^N}\otimes{R}\right),
    \\
    \left(\mu\otimes\mu\right)R=R\left(\mu\otimes\mu\right),
    \\
    \sum_{m=0}^{N-1}\left(R^{\pm1}({id}_{\C^N}\otimes\mu)\right)_{lm}^{km}
    =\left(-q^{1/2}\right)^{\pm1}\delta_{k,l}
  \end{gather*}
with ${id}_{\C^N}$ the identity on $\C^N$.
\end{prop}
\par
Now I will describe how to calculate Kashaev's invariant for knots using this
enhanced Yang-Baxter operator.
\par
\begin{figure}[h]
\includegraphics[scale=0.25]{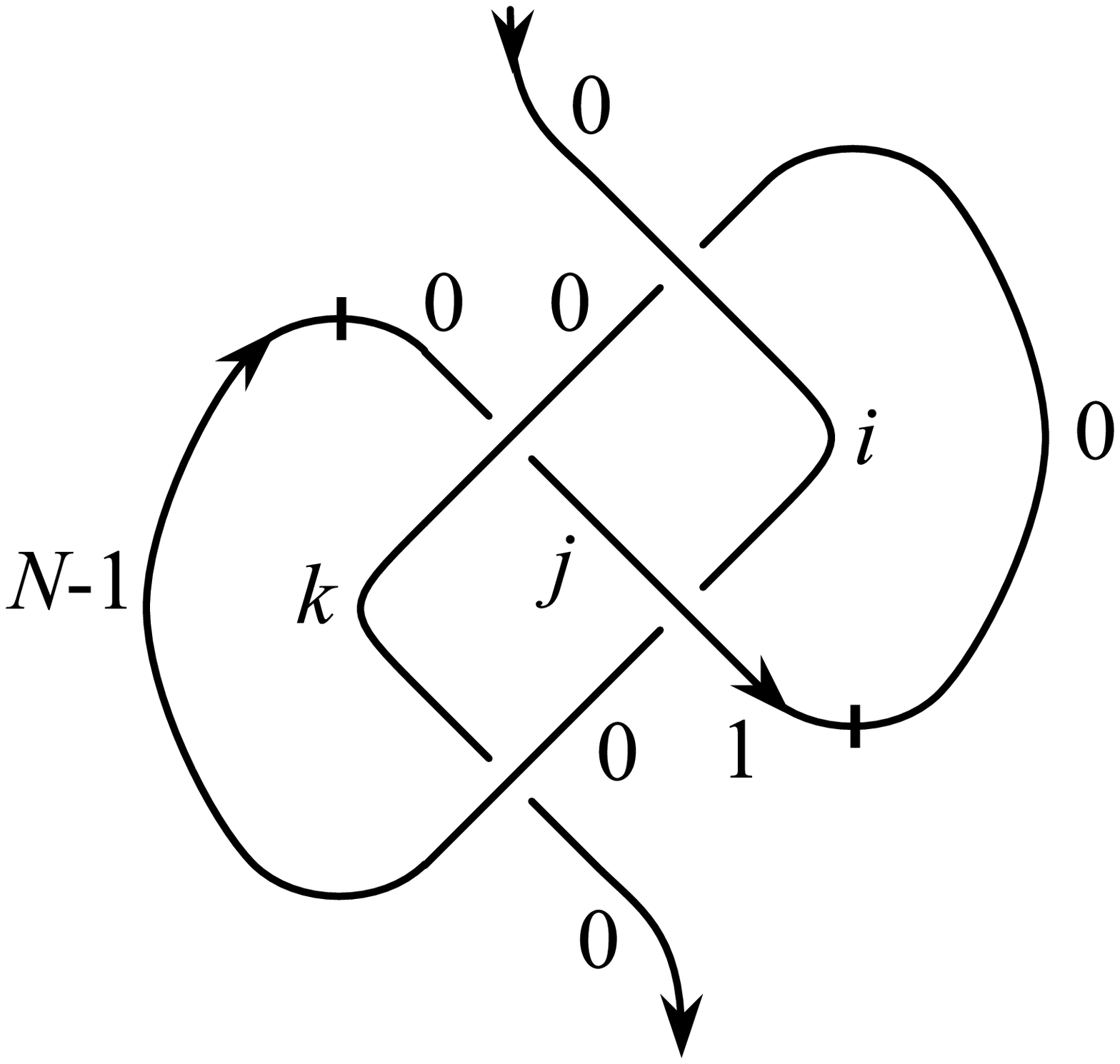}
\caption{Figure-eight knot with a labeling.}
\label{fig:label}
\end{figure}

First we present a knot in a $(1,1)$-tangle formula so that the string comes
from above and goes down and at each crossing both arcs also go down.
We decompose the string into edges so that at each crossing four edges meet and
at each maximum/minimum where the string go from the left to the right two edges
meet.
See Figure~\ref{fig:label} for example.
A {\em labeling} is an assignment of an element in $\{0,1,2,\dots,N-1\}$ to
each edge.
Here we assign $0$ to the edge from the starting point and to the edge to the
end point.
We call the number assigned to an edge its label.
Given a labeling we assign the following quantities to crossings and
maxima/minima:
\setlength{\unitlength}{1mm}
\thicklines
\medskip
\begin{equation*}
\raisebox{-5mm}{\begin{picture}(10,10)
  \put(10,10){\vector(-1,-1){10}}
  \put( 6, 4){\vector( 1,-1){ 4}}
  \put( 4, 6){\line(-1, 1){ 4}}
  \put( 0,10){\makebox(0,0)[br]{$k$}}\put(10,10){\makebox(0,0)[bl]{$n$}}
  \put( 0, 0){\makebox(0,0)[tr]{$l$}}\put(10, 0){\makebox(0,0)[tl]{$m$}}
\end{picture}}
\, :\,R_{lm}^{kn},\qquad
\raisebox{-5mm}{\begin{picture}(10,10)
  \put( 0,10){\vector( 1,-1){10}}
  \put(10,10){\line(-1,-1){ 4}}
  \put( 4, 4){\vector(-1,-1){ 4}}
  \put( 0,10){\makebox(0,0)[br]{$k$}}\put(10,10){\makebox(0,0)[bl]{$n$}}
  \put( 0, 0){\makebox(0,0)[tr]{$l$}}\put(10, 0){\makebox(0,0)[tl]{$m$}}
\end{picture}}
\, :\,\left(R^{-1}\right)_{lm}^{kn},\qquad
\raisebox{-2.5mm}{\begin{picture}(10,5)
  \qbezier(0,5)(0,0)(5,0)
  \qbezier(5,0)(10,0)(10,5)
  \put(5,-1){\line(0,1){2}}
  \put(9,5.5){\vector(0,1){0}}
  \put(1,1){\makebox(0,0)[tr]{$k$}}
  \put(9,1){\makebox(0,0)[tl]{$l$}}
\end{picture}}
\, :\,\mu_{l}^{k},\qquad
\raisebox{-2.5mm}{\begin{picture}(10,5)
  \qbezier(0,0)(0,5)(5,5)
  \qbezier(5,5)(10,5)(10,0)
  \put(5,4){\line(0,1){2}}
  \put(10.1,-0.5){\vector(0,-1){0}}
  \put(1,4){\makebox(0,0)[br]{$k$}}
  \put(9,4){\makebox(0,0)[bl]{$l$}}
\end{picture}}
\,:\,\left(\mu^{-1}\right)_{l}^{k}.
\end{equation*}
\smallskip
\par
The {\em weight} for a labeling is the product of all the
quantities above.
Then the Kashaev invariant $\langle{K}\rangle_N$ of a knot $K$ is defined
to be the sum of all the weights up to some power of $q$ (which does not matter
in our case), where the summation runs over all the labelings.
\par
If four edges meeting at a crossing are labeled $k$, $l$, $m$ and $n$ as above,
we call $[n-m]$, $[k-n]$, $[l-k]$ and $[m-l-1]$ the angles between the adjacent
edges.
Then we have the following proposition.
\begin{prop}[{\cite[Proposition 2]{Yokota:volume}}]\label{prop:angle}
The weight for a labeling vanishes unless all of the following conditions
are satisfied.
For a proof see \cite{Yokota:volume}.
\par
  \begin{enumerate}
  \item the sum of the angles around any crossing is $N-1$.
  \item the sum of the angles contained in any bounded region is $N-1$.
  \item any angle in any of the two unbounded regions is zero.
  \end{enumerate}
\end{prop}
\par
We will calculate the Kashaev invariant for the figure-eight knot.
\par
From Proposition~\ref{prop:angle} (3), we can assume that labels are as in
Figure~\ref{fig:label}.
Then the Kashaev invariant $\langle{4_1}\rangle_{N}$ is
\begin{align*}
  &\sum_{j=0}^{N-1}\sum_{k\in[0,j-1]}\sum_{i\in[1,j]}
  \left(R^{-1}\right)_{0,i}^{0,0}\left(R^{-1}\right)_{0,1}^{j,i}
  R_{N-1,0}^{\phantom{N}k\phantom{1},0}R_{k,j}^{0,0}
  \left(\mu^{-1}\right)_{N-1}^{0}
  \mu_{1}^{0}
  \\
  &=
  N^4\sum_{j=0}^{N-1}\sum_{k\in[0,j-1]}\sum_{i\in[1,j]}
  \frac{q^{-1}}{(q^{-1})_{[i-1]}(q)_{[-i]}}\times
  \frac{q^{-1+(-j)(-i)}}{(q)_{[i-1]}(q^{-1})_{[j-i]}(q)_{[-j]}}
  \\
  &\phantom{N^4\sum_{j=0}^{N-1}\sum_{k\in[0,j-1]}\sum_{i\in[1,j]}}\times
  \frac{q^{1-(N)(-k)}}{(q)_{[k]}(q^{-1})_{[-k-1]}}\times
  \frac{q^{1-(k+1)j}}{(q)_{[j-k-1]}(q^{-1})_{[-j]}(q^{-1})_{[k]}}
  \\
  &=
  N^2\sum_{j=0}^{N-1}
  \frac{q^{-j}}{(q)_{[-j]}(q^{-1})_{[-j]}}
  \sum_{k\in[0,j-1]}
  \frac{q^{-kj}}
       {(q)_{[j-k-1]}(q^{-1})_{[k]}}
  \sum_{i\in[1,j]}
  \frac{q^{ij}}
       {(q)_{[i-1]}(q^{-1})_{[j-i]}}
\end{align*}
up to some power of $q$, since
$(q^{-1})_{[i-1]}(q)_{[-i]}=(q^{-1})_{i-1}(q)_{N-i}
=\dfrac{N}{(q)_{N-i}}(q)_{N-i}=N$
and
$(q)_{[k]}(q^{-1})_{[-k-1]}=(q)_{k}\dfrac{N}{(q)_{k}}=N$.
Here $a\in[b,c]$ means that $q^{b}$, $q^{a}$, and $q^{c}$ are on the unit
circle counterclockwise (any pair may be the same).
Note that $a\in[b,c]$ if and only if $[a-b]+[c-a]=[c-b]$.
\par
From the lemma below we have
\begin{align*}
  \langle{4_1}\rangle_{N}
  &=
  N^2\sum_{j=0}^{N-1}
  \frac{(-1)^{[j-1]}q^{-j-([j-1]+1)([j-1]-2j)/2}}{(q)_{[-j]}(q^{-1})_{[-j]}}
  \sum_{k\in[0,j-1]}
  \frac{q^{-kj}}
       {(q)_{[j-k-1]}(q^{-1})_{[k]}}
  \\
  &=
  N^2\sum_{j=0}^{N-1}
  \frac{1}{(q)_{[-j]}(q^{-1})_{[-j]}}
  \\
  &=
  \sum_{j=0}^{N-1}(q)_{[N-1+j]}(q^{-1})_{[N-1+j]}
  \\
  &=
  \sum_{j=0}^{N-1}(q)_{j}(q^{-1})_{j}.
\end{align*}
Note that the final formula coincides with Kashaev's
\cite[(2.2)]{Kashaev:LETMP97}.
\par
The following lemma is a key to reduce the summations of the Kashaev's
invariant.
\begin{lem}[Yokota's version of {\cite[Lemma~3.2]{Murakami/Murakami:volume}}]
Let $0\le{k},{l},{m}\le{N-1}$ be integers.
Then
\begin{equation*}
  \sum_{k\in[l,m]}
  \frac{q^{-(m-l+1)k}}{(q)_{[m-k]}(q^{-1})_{[k-l]}}
  =
  (-1)^{[m-l]}q^{([m-l]+1)([m-l]-2m)/2}.
\end{equation*}
\end{lem}
\begin{proof}
We use the following equality \cite[Lemma~3.2]{Murakami/Murakami:volume}.
\begin{equation*}
\begin{split}
  &\sum_{i=0}^{N-1-\alpha}
  \frac{q^{(\beta-\alpha)i/2}}{(q)_{i}(q^{-1})_{N-1-\alpha-i}}
  \\
  &\quad=\frac{(-1)^{N-1-\alpha}q^{-N(N-1)/2-\beta(\alpha+1)/2}}{N}
   \frac{(q)_{\alpha}(q)_{N-1-\alpha+[(\alpha-\beta)/2]}}
        {(q)_{[(\alpha-\beta)/2]}}.
\end{split}
\end{equation*}
\par
To prove the lemma, we put $\beta=-\alpha$.
Then we have
\begin{equation*}
  \sum_{i=0}^{N-1-\alpha}
  \frac{q^{-\alpha i}}{(q)_{i}(q^{-1})_{N-\alpha-i-1}}
  =
  (-1)^{N-1-\alpha}q^{-N(N-1)/2+\alpha(\alpha+1)/2}.
\end{equation*}
Putting $i=[m-k]$ and $N-\alpha-1=[m-l]$, the conclusion follows.
\end{proof}
\section{Geometry}\label{sec:geometry}
In this section I will describe how to decompose a knot complement
$S^3\setminus{K}$ into ideal topological tetrahedra and how these
define hyperbolic structure when $K$ is the figure-eight knot.
Here `ideal' means that each vertex is at the cusp (the infinite point
corresponding to $K$).
(That is, all vertices of the tetrahedra meet at a point, and the complement of
the regular neighborhood of the point is $S^3\setminus{K}$.)
\par
Fix a knot diagram throughout.
\subsection{Ideal tetrahedron decomposition of
$S^3\setminus({K\cup\text{two points}})$}
We will decompose $S^3\setminus({K\cup\text{two points}})$ into tetrahedra
so that all the vertices are either on $K$ or on the two points outside $K$.
\par
\begin{figure}[h]
\includegraphics[scale=0.45]{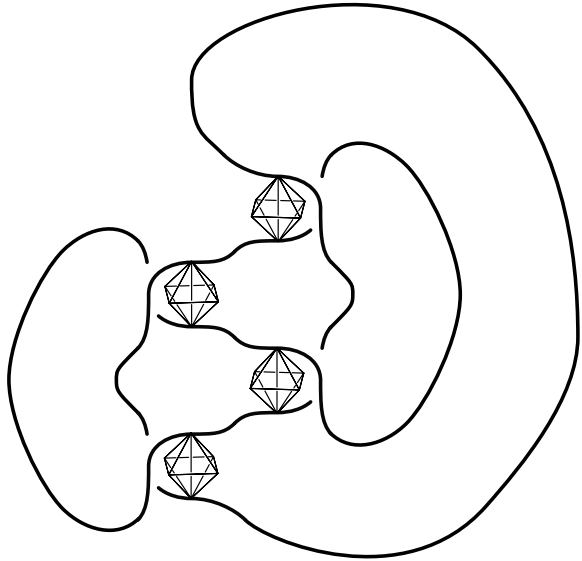}
\caption{Put octahedra at crossings.}
\label{fig:fig8_octa}
\end{figure}
\begin{figure}[h]
\includegraphics[scale=0.3]{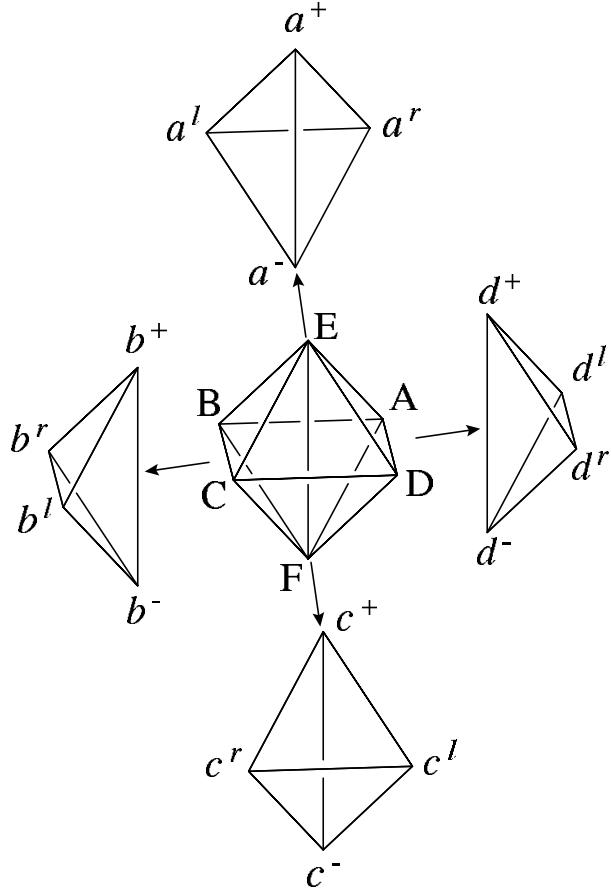}
\caption{Decompose an octahedron into four tetrahedra.}
\label{fig:octa}
\end{figure}
\begin{figure}[h]
\includegraphics[scale=0.3]{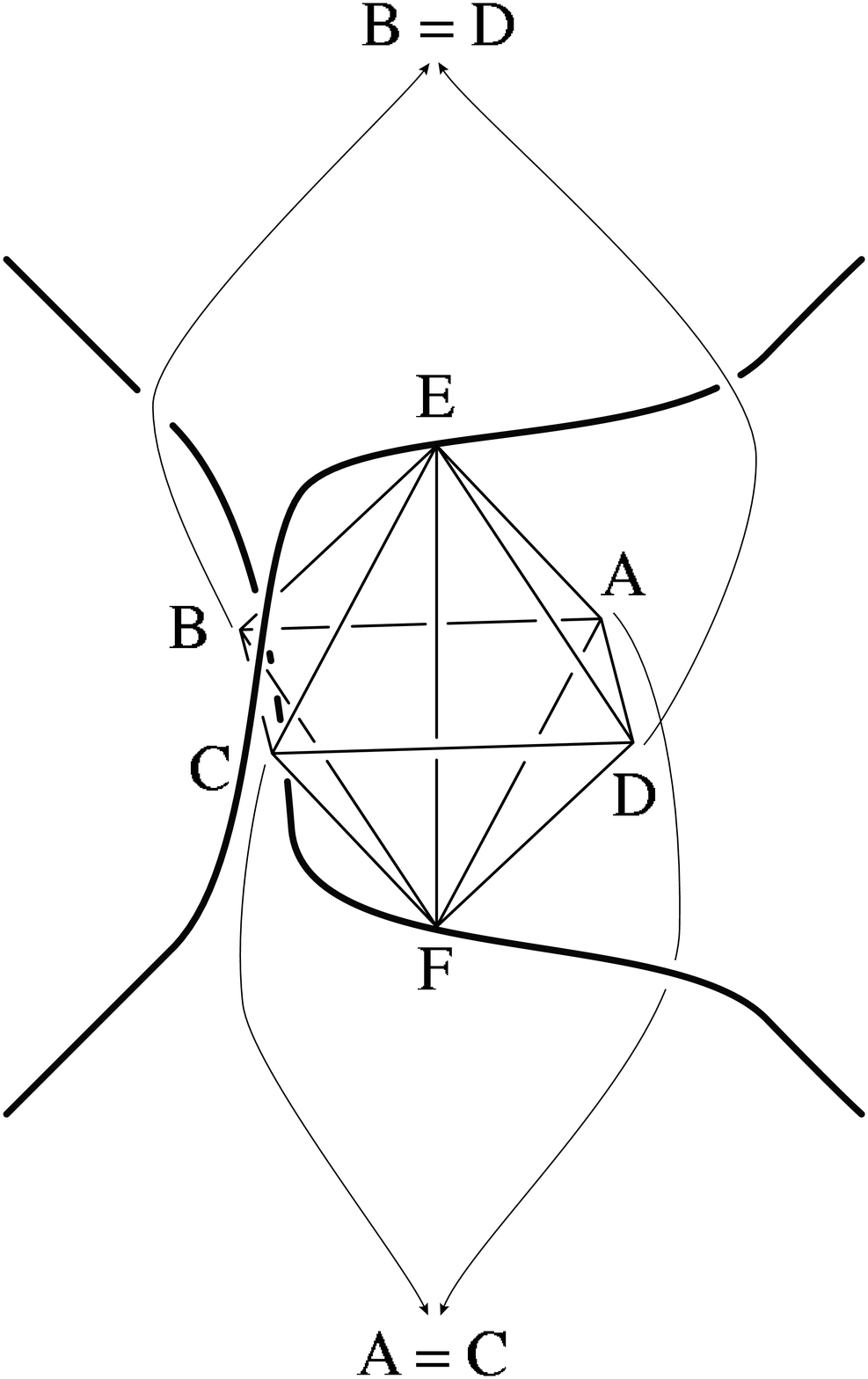}
\caption{Pull the pairs of vertices to $\pm\infty$.}
\label{fig:octa_plus_minus}
\end{figure}

First we put an octahedron
${\rm E}_i$-${\rm A}_i{\rm B}_i{\rm C}_i{\rm D}_i$-${\rm F}_i$
at the $i$th crossing of the knot diagram so that ${\rm E}_i{\rm F}_i$ is
perpendicular to the plane, that ${\rm E}_i$ and ${\rm F}_i$ touch the
over-crossing point and the under-crossing point respectively,
and that ${\rm A}_i$, ${\rm B}_i$, ${\rm C}_i$, and ${\rm D}_i$ are to
the northeast, northwest, southwest, and southeast, respectively.
(see Figures~\ref{fig:fig8_octa} and \ref{fig:octa_plus_minus}).
We call this the $i$th octahedron.
We decompose each octahedron into four tetrahedra as shown in
Figure~\ref{fig:octa}.
(For the figure-eight knot, we now have 16 tetrahedra.)
We name each vertex of the tetrahedra as in Figure~\ref{fig:octa}.
(Here we drop subscriptions.)
\par
\begin{figure}[h]
\includegraphics[scale=0.25]{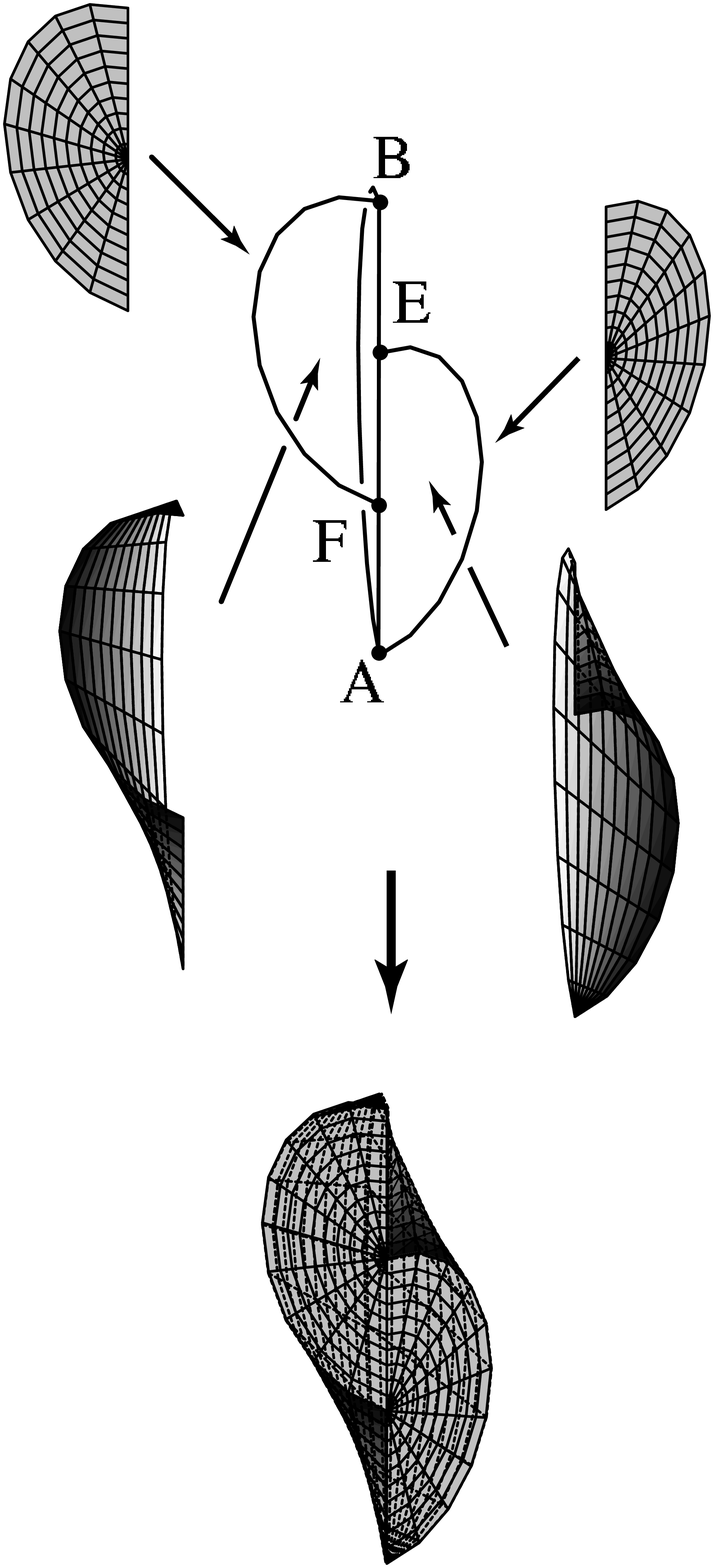}\qquad\qquad
\includegraphics[scale=0.25]{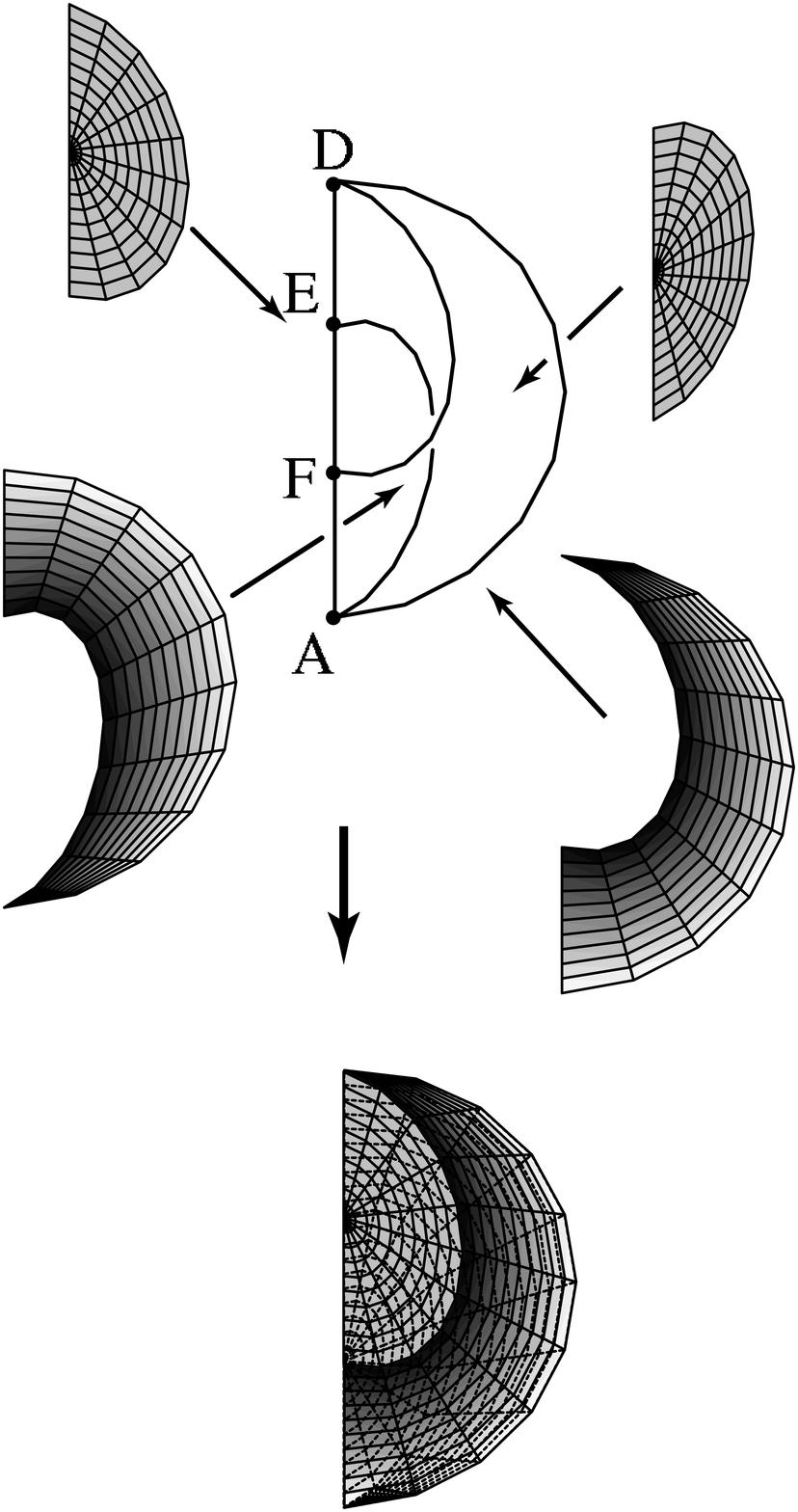}
\includegraphics[scale=0.25]{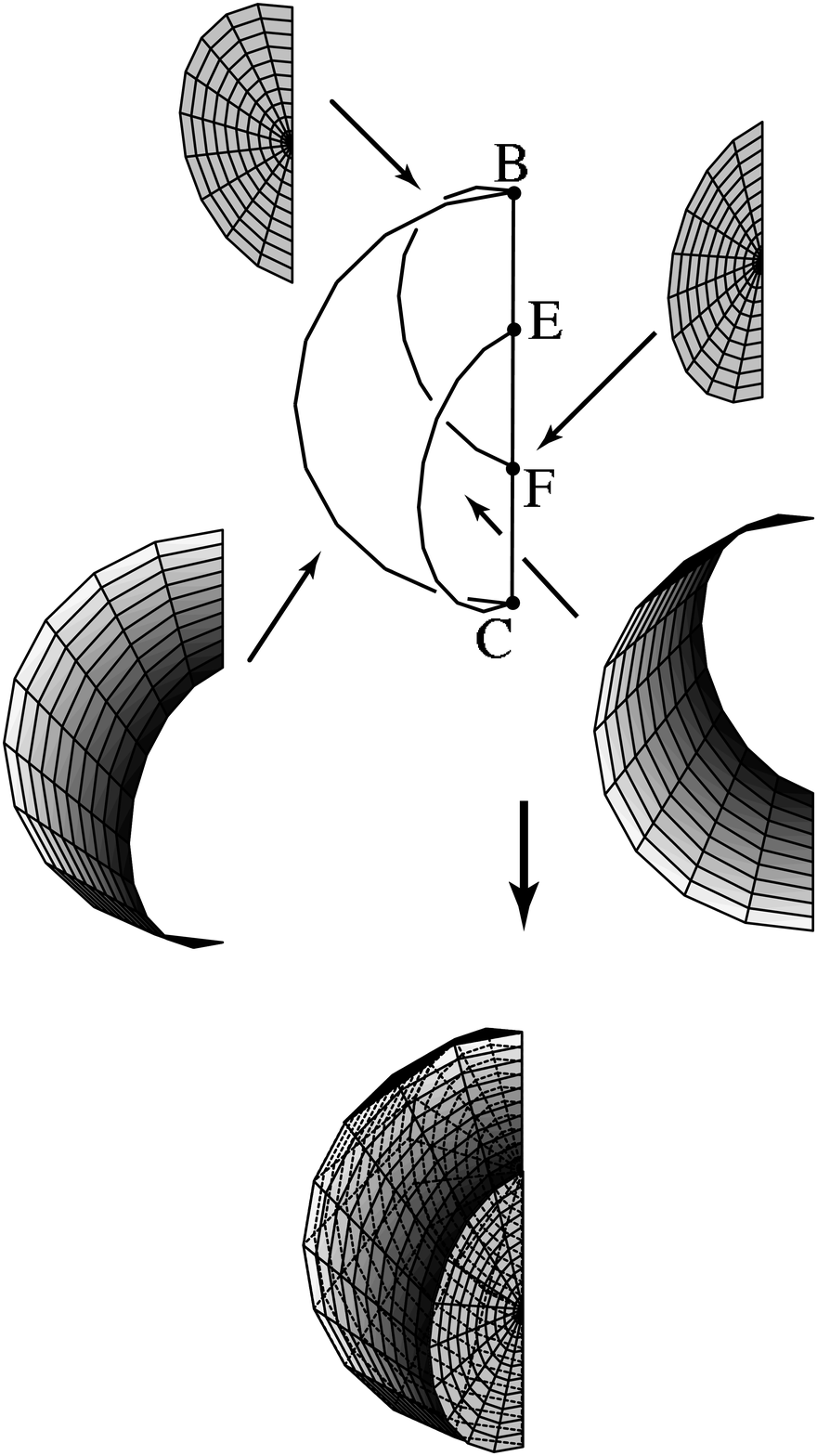}\qquad\qquad
\includegraphics[scale=0.25]{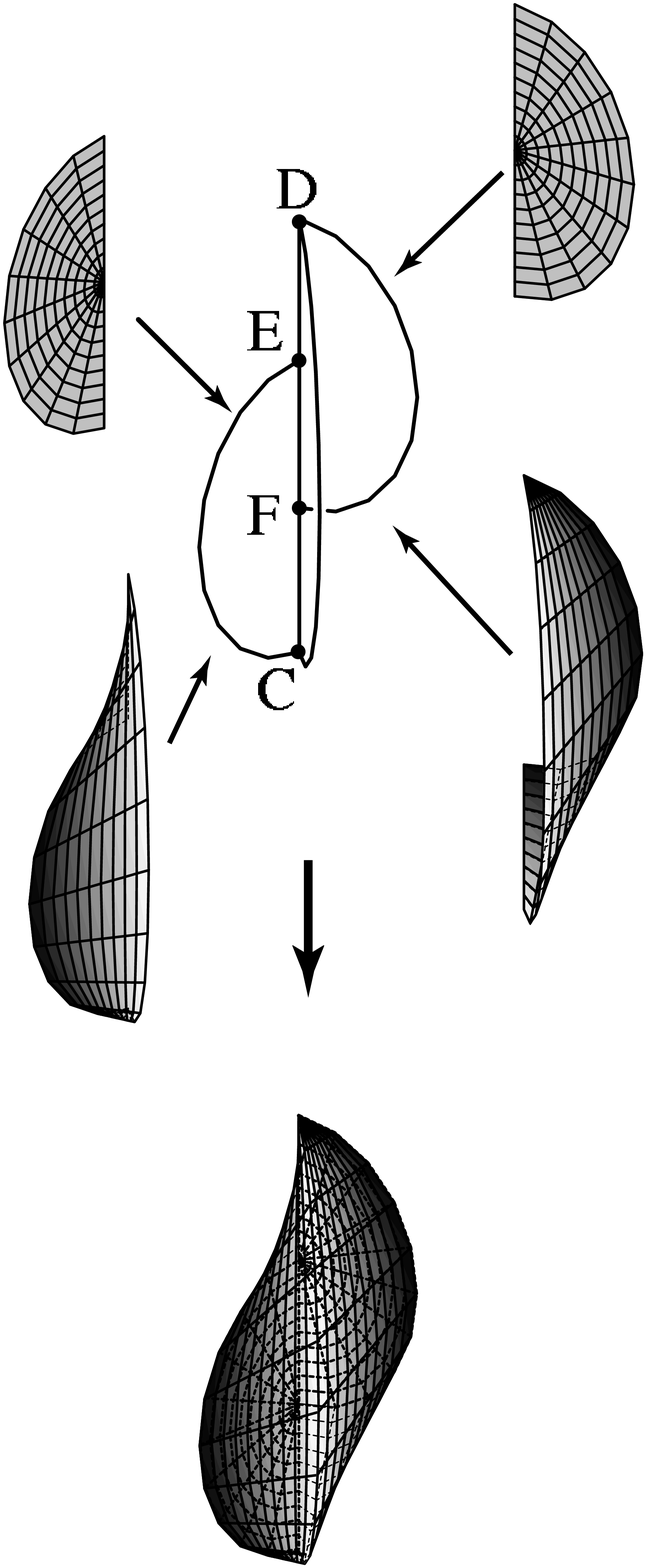}
\caption{Twisted tetrahedra.}
\label{fig:twisted_tetra}
\end{figure}
\begin{figure}[h]
\includegraphics[scale=0.35]{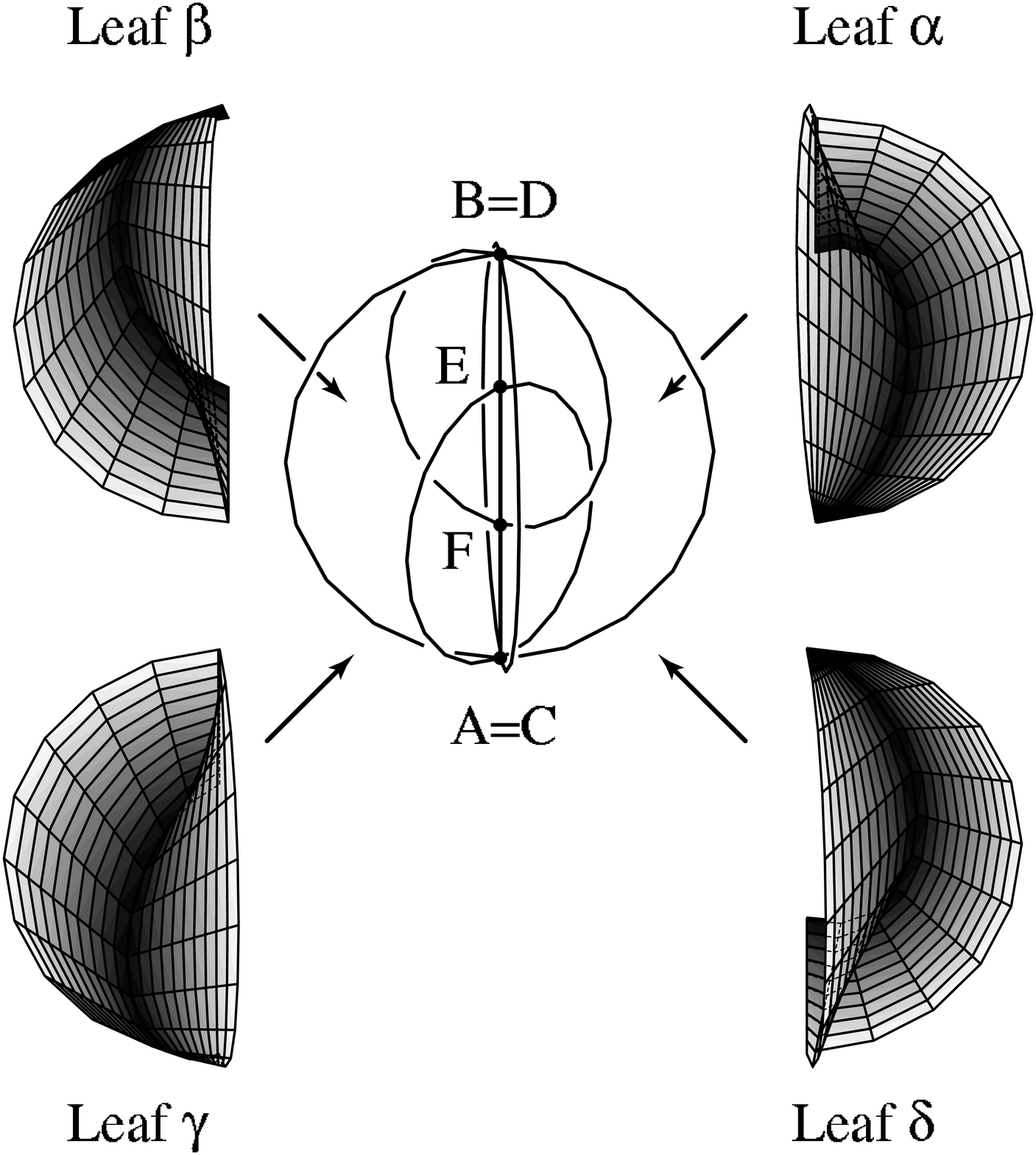}
\caption{A twisted octahedron and its four leaves.}
\label{fig:twisted_octa}
\end{figure}

We deform an octahedron attached to a positive crossing as follows.
We pull the vertices ${\rm B}$ and ${\rm D}$ upward and identify
the edge ${\rm EB}$ with the edge ${\rm ED}$.
Similarly we pull the vertices ${\rm A}$ and ${\rm C}$ downward and identify
the edge ${\rm FA}$ with the edge ${\rm FC}$.
(See Figure~\ref{fig:octa_plus_minus}.)
Now the resulting {\em twisted octahedron} consists of four
{\em twisted tetrahedra} (Figure~\ref{fig:twisted_tetra}).
The surface of a twisted octahedron consists of four {\em leaves} made from
$\triangle{b^+}{b^r}{b^l}\cup\triangle{c^+}{c^r}{c^l}$,
$\triangle{c^-}{c^l}{c^r}\cup\triangle{d^-}{d^l}{d^r}$,
$\triangle{d^+}{d^r}{d^l}\cup\triangle{a^+}{a^r}{a^l}$,
and
$\triangle{a^-}{a^l}{a^r}\cup\triangle{b^-}{b^l}{b^r}$,
which we will call Leaf~$\gamma$, Leaf~$\delta$, Leaf~$\alpha$, and
Leaf~$\beta$, respectively (Figure~\ref{fig:twisted_octa}).
\par
\begin{figure}[h]
\includegraphics[scale=0.35]{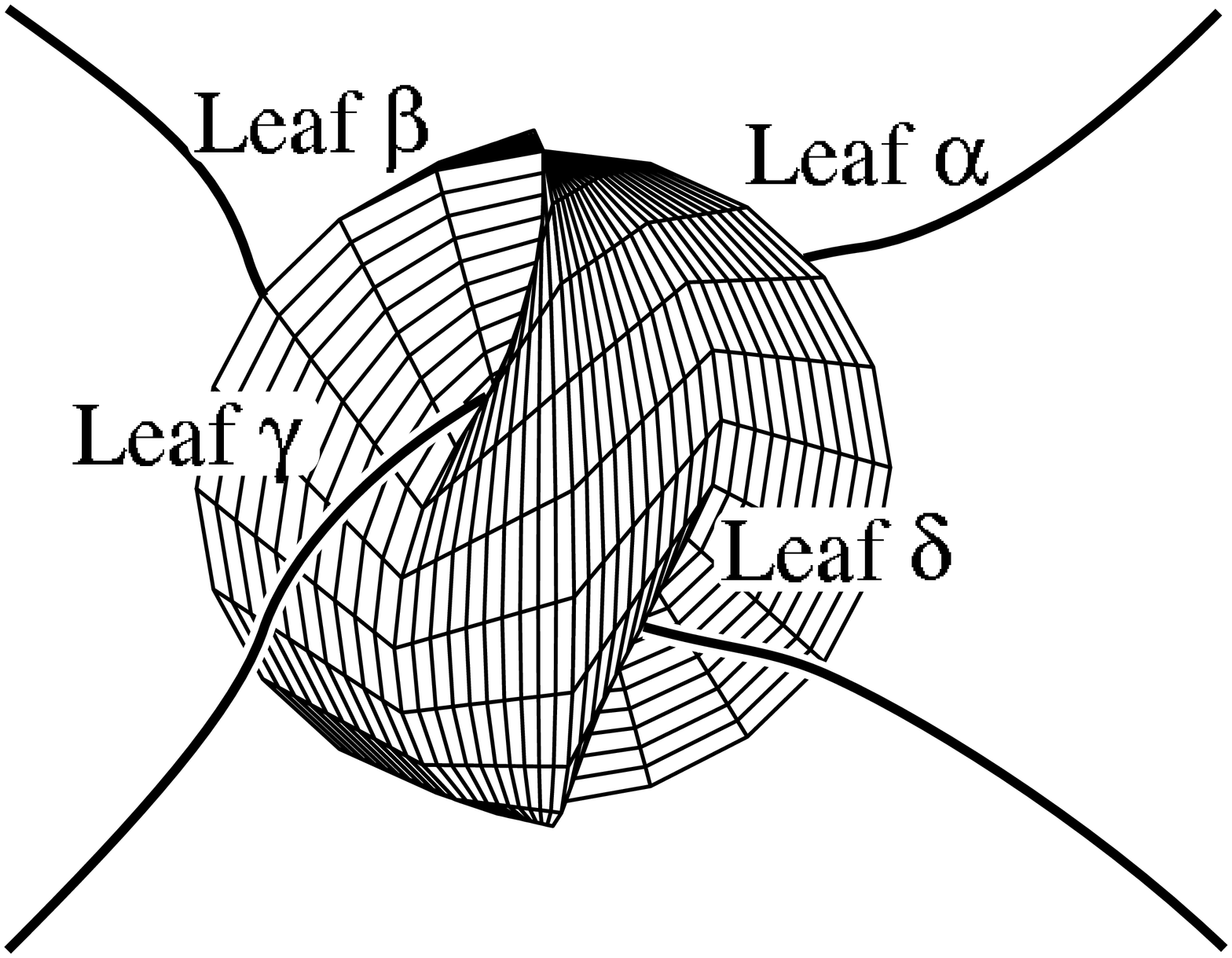}
\caption{A twisted octahedron and a crossing.}
\label{fig:twisted_octa_cross}
\end{figure}

Now this twisted octahedron and the near positive crossing look like
Figure~\ref{fig:twisted_octa_cross}.
\par
For an octahedron attached to a negative crossing we deform it similarly.
We pull the vertices ${\rm B}$ and ${\rm D}$ downward, and ${\rm A}$ and
${\rm C}$ upward to identify the edge ${\rm EB}$ with the edge ${\rm ED}$,
and the edge ${\rm FA}$ with the edge ${\rm FC}$.
In this case Leaf~$\alpha$ is made from $\triangle{d^-}{d^l}{d^r}$ and
$\triangle{a^-}{a^l}{a^r}$,
and Leaves~$\beta$, $\gamma$, and $\delta$ are from
$\triangle{a^+}{a^r}{a^l}$ and $\triangle{b^+}{b^r}{b^l}$,
$\triangle{b^-}{b^l}{b^r}$ and $\triangle{c^-}{c^l}{c^r}$, and
$\triangle{c^+}{c^r}{c^l}$ and $\triangle{d^+}{d^r}{d^l}$, respectively.
\par
\begin{figure}[h]
\includegraphics[scale=0.34]{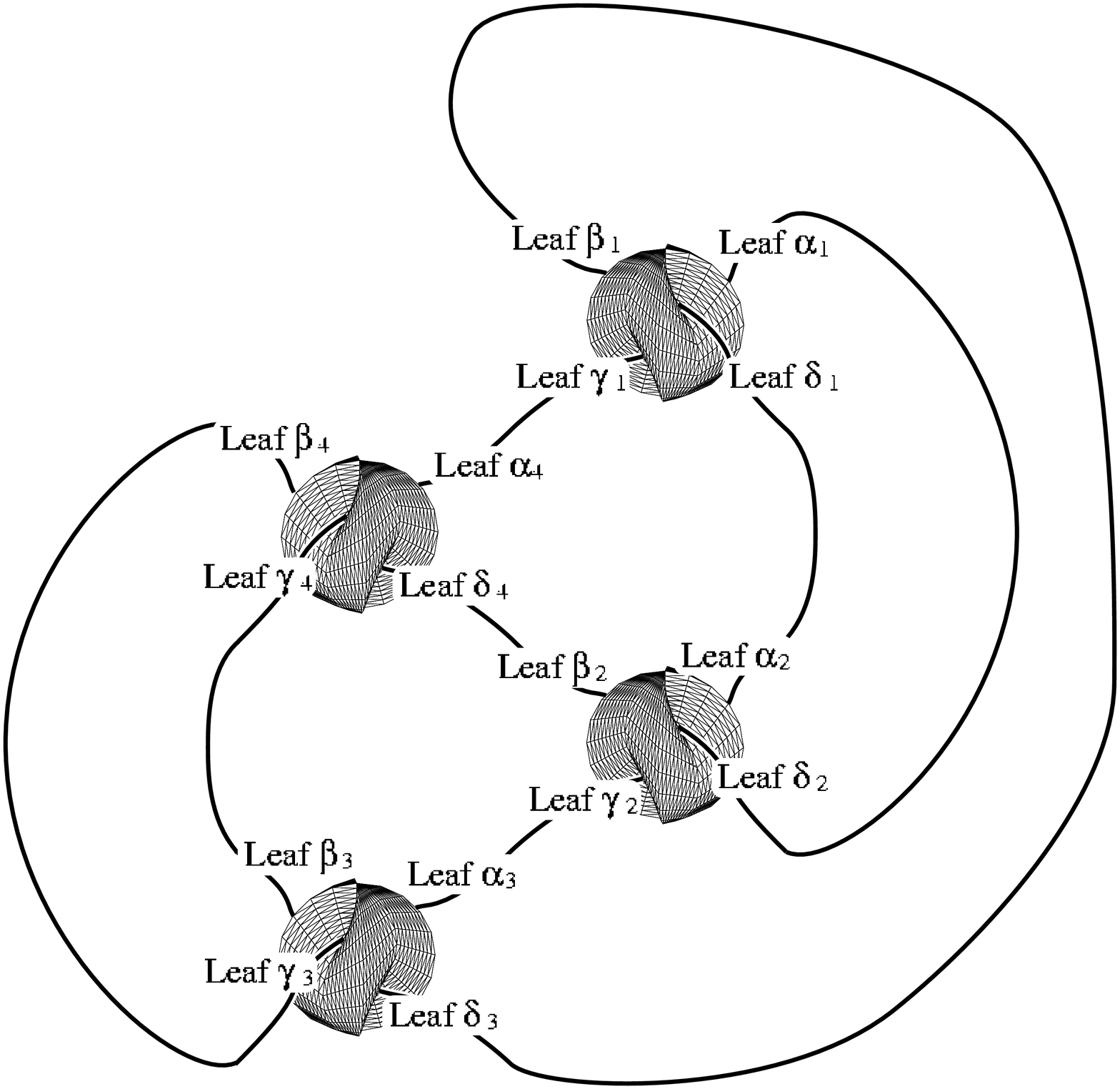}
\caption{Figure-eight knot with twisted octahedra.}
\label{fig:fig8_twisted_octa}
\end{figure}

Now the whole picture for the figure-eight knot is as shown in
Figure~\ref{fig:fig8_twisted_octa}.
\par
Next we identify two leaves facing each other along an arc by thickening
the octahedra and shortening the arc.
For the figure-eight knot case, we identify Leaves
$\delta_1$ and $\alpha_2$,
$\gamma_2$ and $\alpha_3$,
$\gamma_3$ and $\beta_4$,
$\delta_4$ and $\beta_2$,
$\delta_2$ and $\alpha_1$,
$\gamma_1$ and $\alpha_4$,
$\gamma_4$ and $\beta_3$, and
$\delta_3$ and $\beta_1$.
Then the eight vertices ${\rm A_1}$, ${\rm C_1}$, ${\rm A_2}$, ${\rm C_2}$,
${\rm B_3}$, ${\rm D_3}$, ${\rm B_4}$, and ${\rm D_4}$ are identified.
We call the resulting vertex $+\infty$.
Similarly the eight vertices ${\rm B_1}$, ${\rm D_1}$, ${\rm B_2}$, ${\rm D_2}$,
${\rm A_3}$, ${\rm C_3}$, ${\rm A_4}$, and ${\rm C_4}$ are identified with the
new vertex $-\infty$.
Moreover the other eight vertices $E_1, F_1,\dots,E_4,F_4$ are also identified,
which we call ${\rm G}$.
\par
\begin{figure}[h]
\includegraphics[scale=0.44]{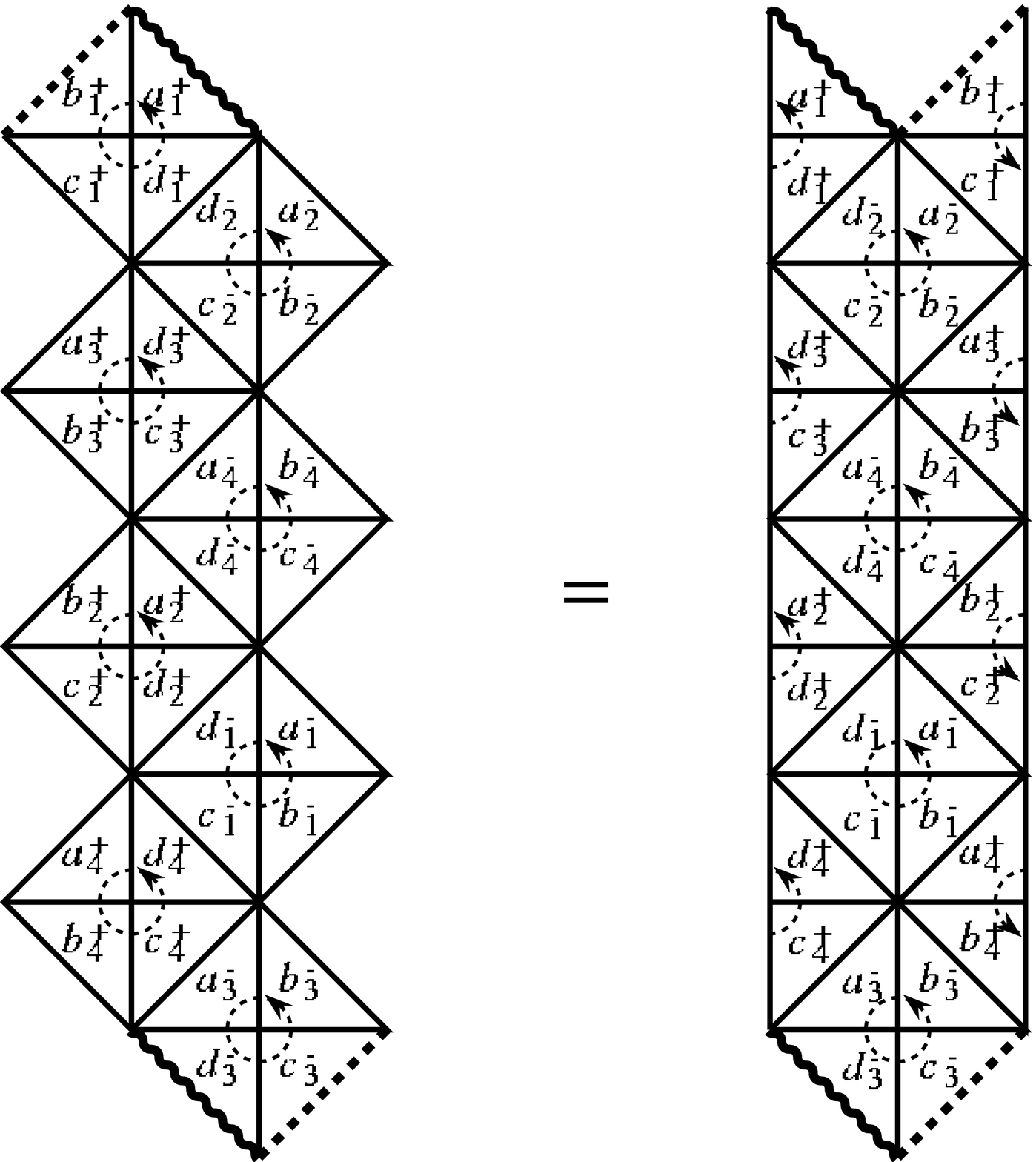}
\caption{Boundary of the regular neighborhood of ${\rm G}$.
(The wavy lines, the dotted lines, and the sidelines are identified
respectively.)}
\label{fig:nbd_K}
\end{figure}
\begin{figure}[h]
\includegraphics[scale=0.33]{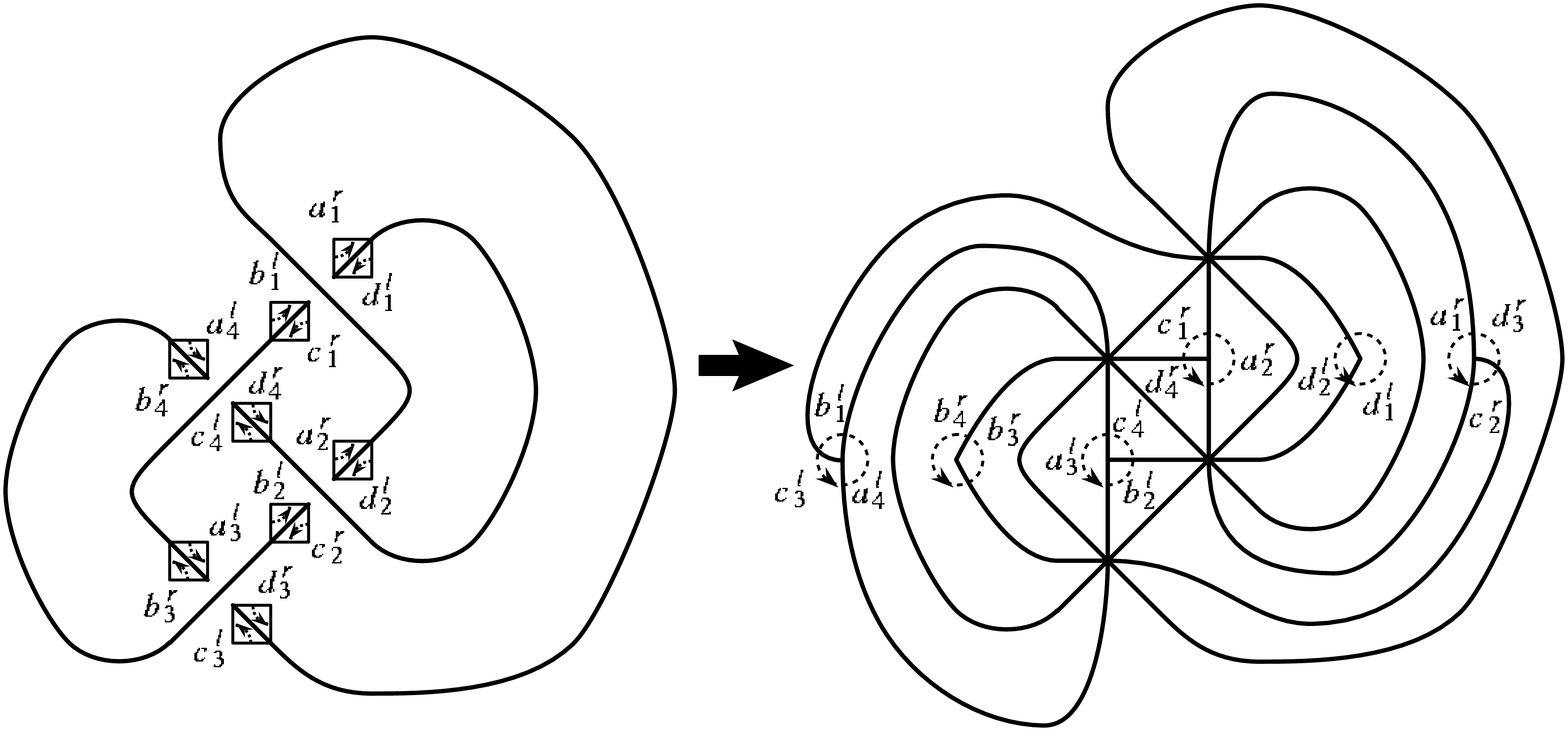}
\caption{Boundaries of the regular neighborhoods of
${\rm A_1=C_1}$, ${\rm A_2=C_2}$, ${\rm B_3=D_3}$, and ${\rm B_4=D_4}$
together with knot diagram before thickening the twisted octahedra viewed from
the top (left).
Boundary of the regular neighborhood of $+\infty$ (right).}
\label{fig:nbd_plus}
\end{figure}
\begin{figure}[h]
\includegraphics[scale=0.33]{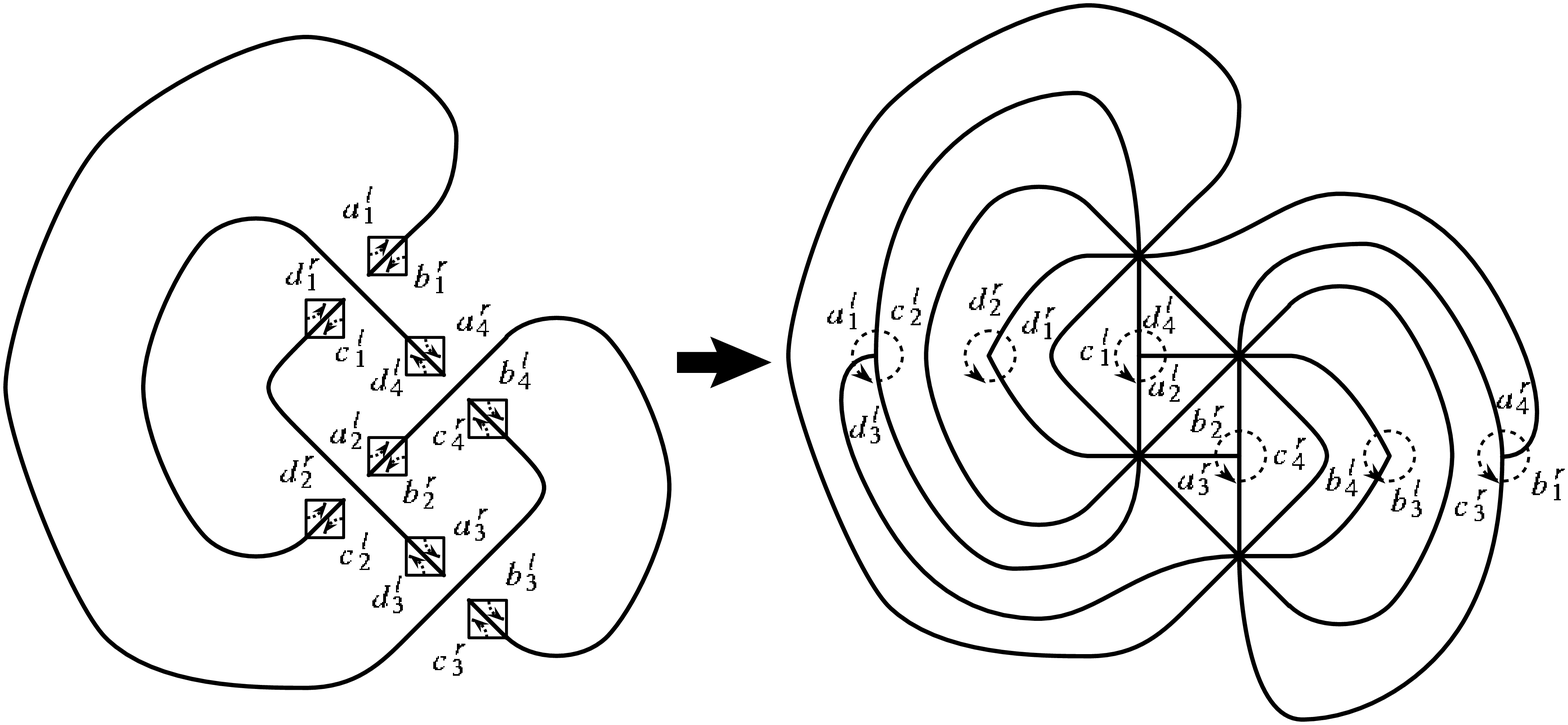}
\caption{Boundaries of the regular neighborhoods of
${\rm B_1=D_1}$, ${\rm B_2=D_2}$, ${\rm A_3=C_3}$, and ${\rm A_4=C_4}$
together with knot diagram before thickening the twisted octahedra viewed from
the bottom (left).
Boundary of the regular neighborhood of $-\infty$ (right).}
\label{fig:nbd_minus}
\end{figure}

Now we have an ideal tetrahedron decomposition of
$S^3\setminus({K\cup\pm\infty})$ with three ideal vertices.
The boundary torus of the regular neighborhood of ${\rm G}$ and the boundary
spheres of the regular neighborhoods of $\pm\infty$ look like
Figures~\ref{fig:nbd_K}, \ref{fig:nbd_plus}, and \ref{fig:nbd_minus}
respectively if $K$ is the figure-eight knot.
Note that we always look at boundaries from the outside.
\subsection{Ideal tetrahedron decomposition of $S^3\setminus{K}$}
We will deform the ideal tetrahedron decomposition of
$S^3\setminus({K\cup\pm\infty})$ given above to obtain an ideal tetrahedron
decomposition of $S^3\setminus{K}$.
\par
To do this we collapse Leaf~$\beta_1$ (= Leaf~$\delta_3$) into the point
$G=E_1=F_1=\dots=E_4=F_4$ and extend this collapsing over the other tetrahedra
linearly.
Note that this leaf corresponds to the broken point of the knot diagram to make
it a $(1,1)$-tangle in \S~\ref{sec:algebra}.
\par
I will show how octahedra will collapse.
\par
\begin{figure}[h]
\includegraphics[scale=0.23]{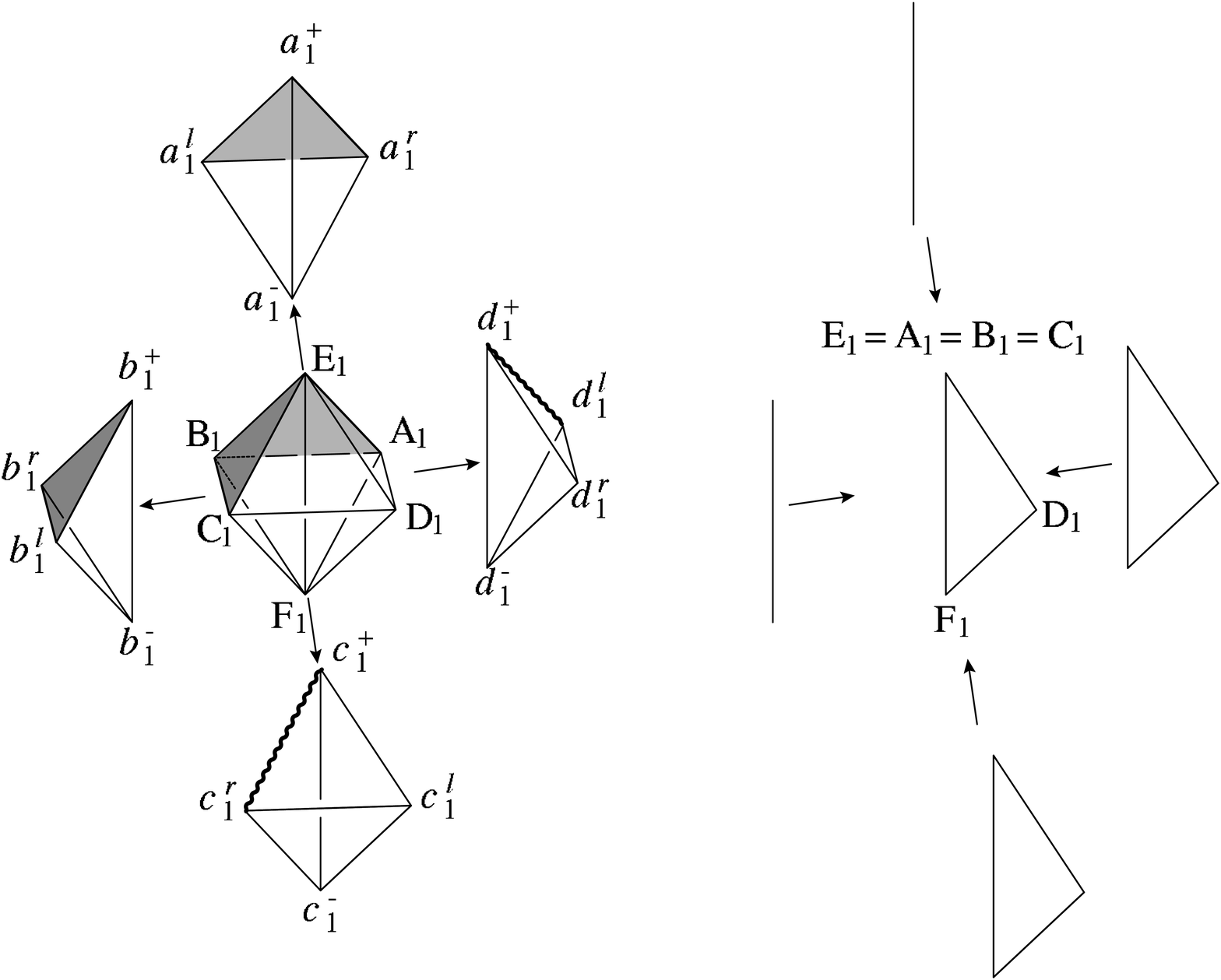}
\caption{The first octahedron is collapsed to the triangle ${\rm F_1E_1D_1}$.
(The shaded triangles and the wavy lines are collapsed.)}
\label{fig:collapsed_octa1}
\end{figure}

In the first octahedron ${\rm E_1}$-${\rm A_1B_1C_1D_1}$-${\rm F_1}$, there are
two faces $\triangle{a^+_1}{a^r_1}{a^l_1}$ and
$\triangle{b^+_1}{b^r_1}{b^l_1}$ which make Leaf~$\beta_1$.
Therefore the tetrahedron ${\rm F_1E_1A_1B_1}$ and ${\rm F_1E_1B_1C_1}$
(they are the cones over $\triangle{a^+_1}{a^r_1}{a^l_1}$ and
$\triangle{b^+_1}{b^r_1}{b^l_1}$ with ${\rm F_1}$)
are collapsed to edges and so three edges ${\rm F_1A_1}$, ${\rm F_1B_1}$ and
${\rm F_1C_1}$ are identified with the edge ${\rm F_1E_1}$.
The tetrahedron ${\rm F_1E_1C_1D_1}$ and ${\rm F_1E_1D_1A_1}$
(they are the joins ${\rm E_1C_1}*{\rm D_1F_1}$ and ${\rm E_1A_1}*{\rm D_1F_1}$)
are collapsed to triangles and so two triangles ${\rm F_1D_1C_1}$ and
${\rm F_1A_1D_1}$ are identified with the triangle ${\rm E_1F_1D_1}$.
So the first octahedron is collapsed to the triangle ${\rm F_1E_1D_1}$.
See Figure~\ref{fig:collapsed_octa1}.
\par
Similarly the third octahedron is collapsed to the triangle ${\rm F_3E_3B_3}$
since it contains the faces $\triangle{c^-_3}{c^l_3}{c^r_3}$ and
$\triangle{d^-_3}{d^l_3}{d^r_3}$ which make
Leaf~$\delta_3$.
Here two triangles ${\rm E_3B_3C_3}$ and ${\rm E_3A_3B_3}$ are identified with
the triangle ${\rm F_3E_3B_3}$.
\par
\begin{figure}[h]
\includegraphics[scale=0.23]{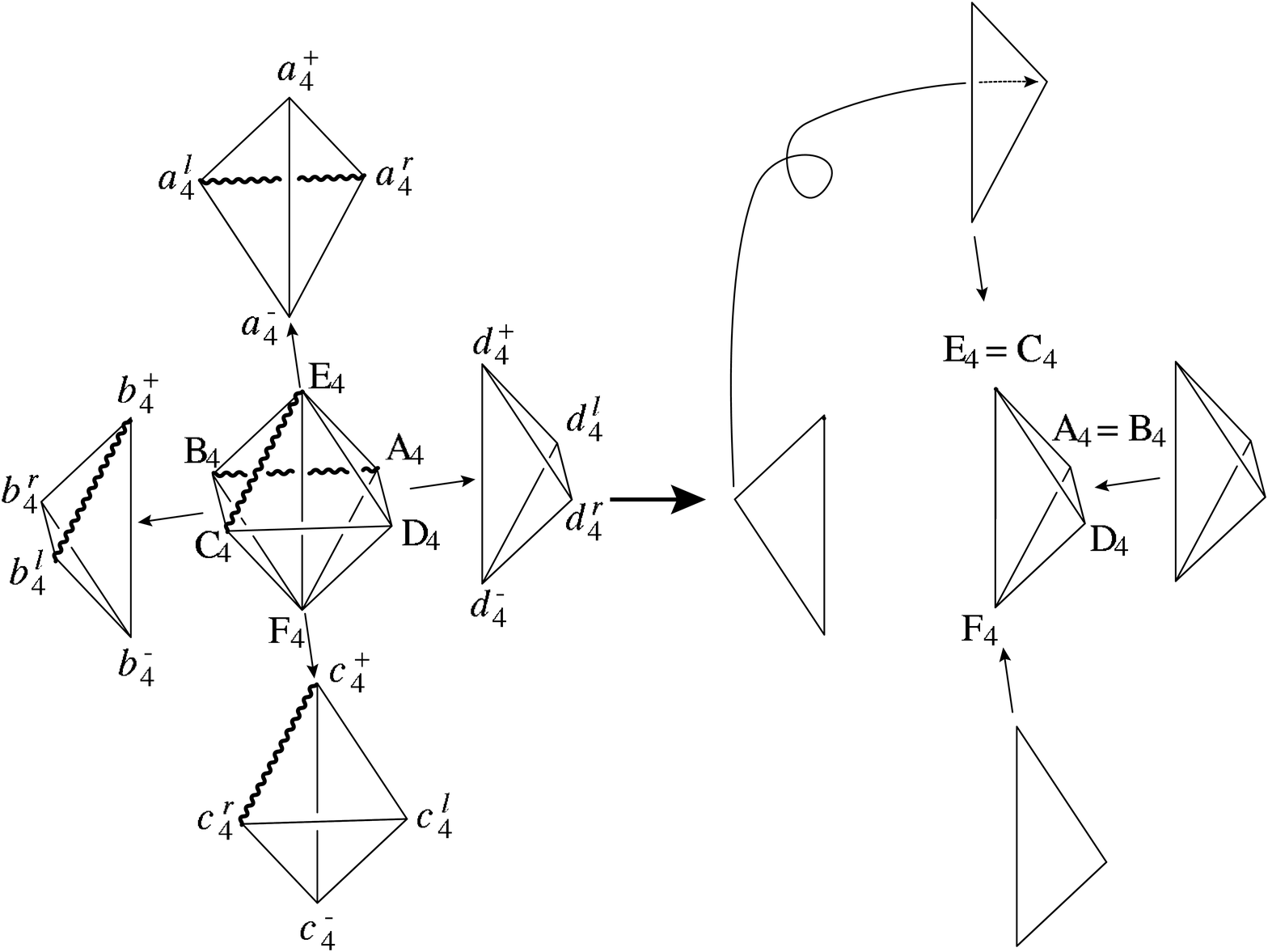}
\caption{The fourth octahedron is collapsed to the tetrahedron
${\rm F_4E_4D_4A_4}$.
(The wavy lines are collapsed.)}
\label{fig:collapsed_octa4}
\end{figure}

On the other hand the edge ${\rm E_4C_4}$ is identified with the edge
${\rm C_3F_3}={\rm A_3F_3}$, which is contained in Leaf~$\delta_3$.
Moreover the edge ${\rm B_4A_4}$ is identified with the edge ${\rm C_1B_1}$,
which is contained in Leaf~$\beta_1$.
Therefore the fourth octahedron is collapsed to the tetrahedron ${\rm F_4E_4D_4A_4}$
as in Figure~\ref{fig:collapsed_octa4}.
Here the tetrahedron ${\rm F_4E_4B_4C_4}$ is collapsed to the triangle
${\rm F_4E_4B_4}$, which is identified with the triangle ${\rm F_4E_4B_4}$.
Then the tetrahedron ${\rm F_4E_4A_4B_4}$ is collapsed to the triangle
${\rm F_4E_4A_4}$.
The tetrahedron ${\rm F_4E_4C_4D_4}$ is collapsed to the triangle
${\rm F_4E_4D_4}$.
\par
Similarly the second octahedron is collapsed to the tetrahedron
${\rm F_2E_2B_2C_2}$ since two edges ${\rm E_2A_2}$ and ${\rm C_2D_2}$ are
collapsed to points respectively.
\par
After all only two tetrahedron ${\rm F_2E_2B_2C_2}$ and ${\rm F_4E_4D_4A_4}$
survive.
\par
\begin{figure}[h]
\includegraphics[scale=0.23]{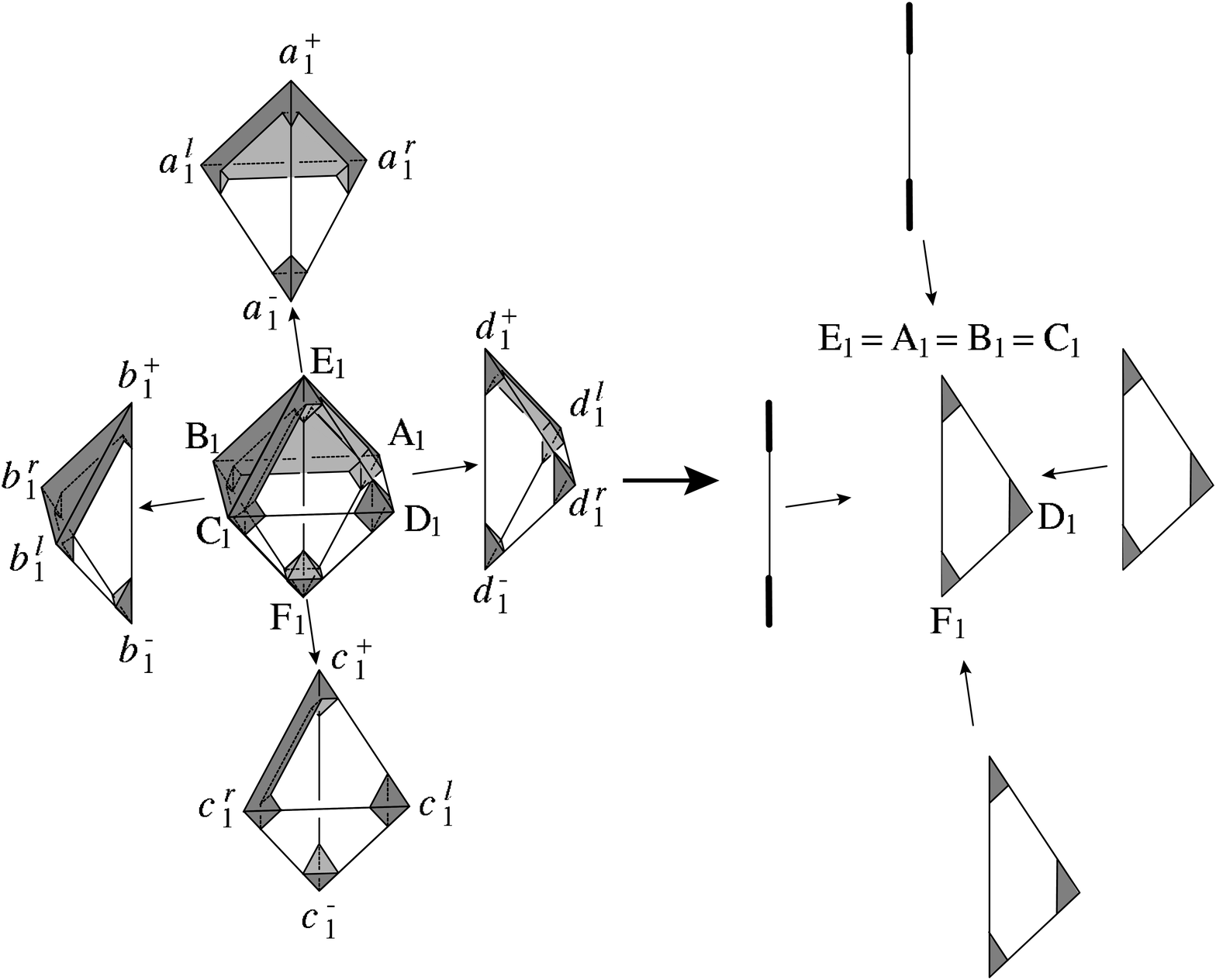}
\caption{The regular neighborhood of Leaf~$\beta_1$ (= Leaf~$\delta_3$)
in the first octahedron.}
\label{fig:nbd_octa1}
\end{figure}
\begin{figure}[h]
\includegraphics[scale=0.23]{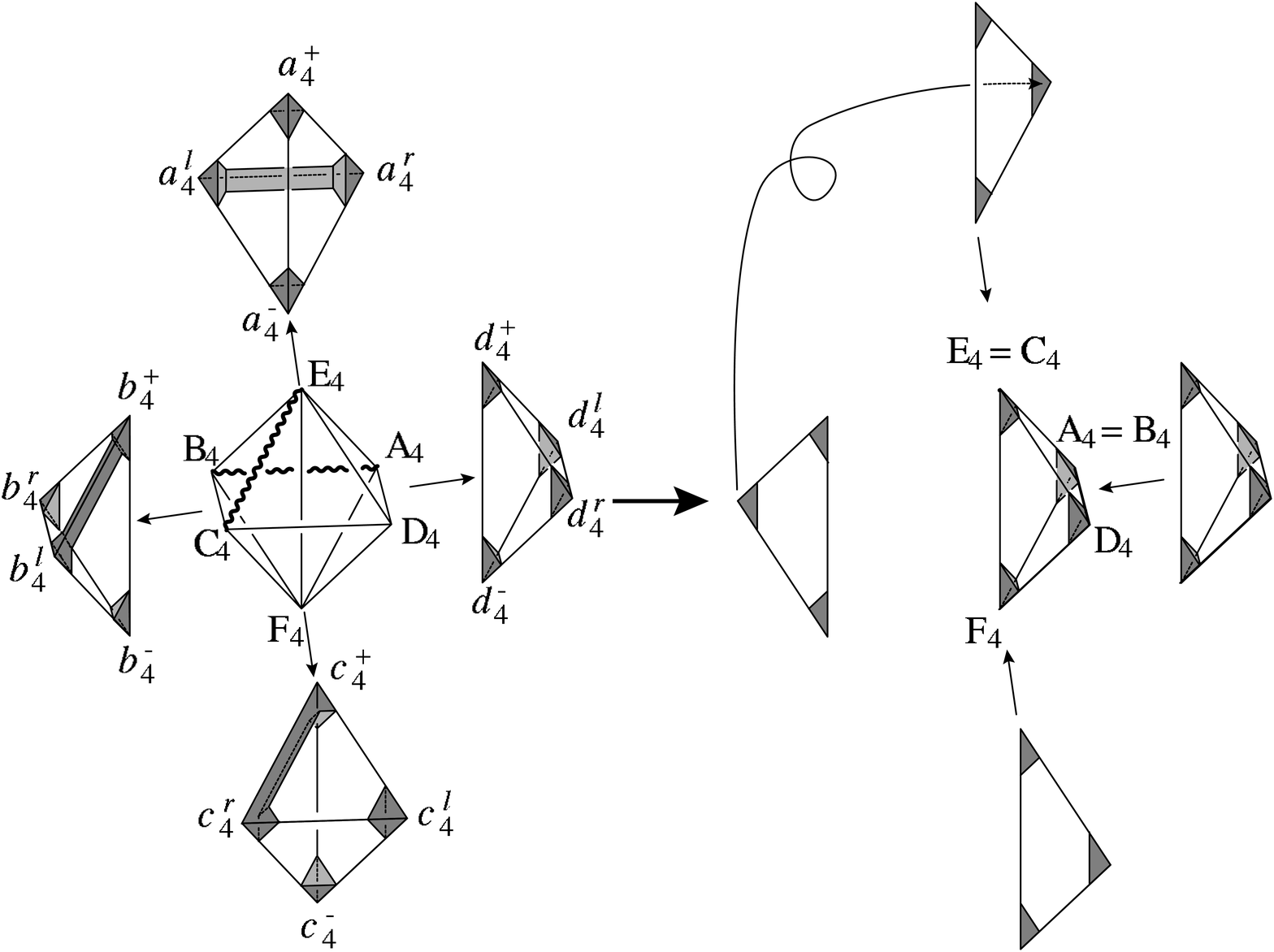}
\caption{The regular neighborhood of Leaf~$\beta_1$ (= Leaf~$\delta_3$)
in the fourth octahedron.}
\label{fig:nbd_octa4}
\end{figure}
\begin{figure}[h]
\includegraphics[scale=0.45]{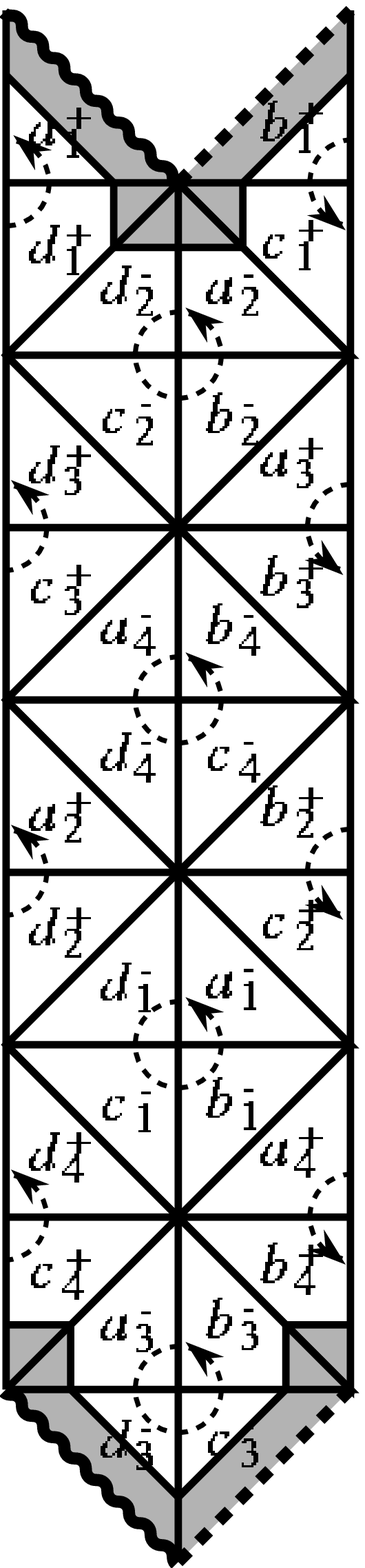}\\
\includegraphics[scale=0.3]{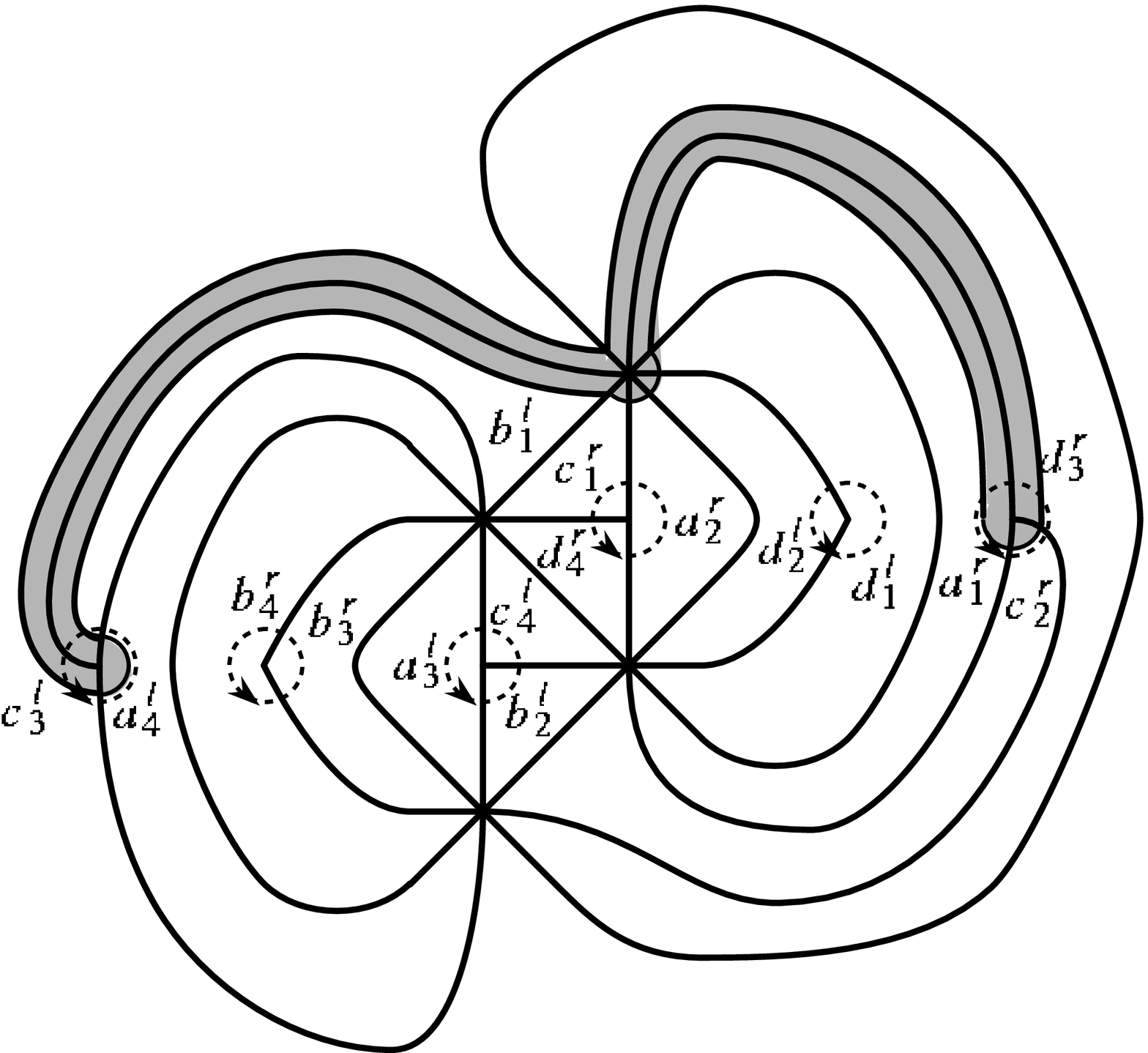}\qquad
\includegraphics[scale=0.3]{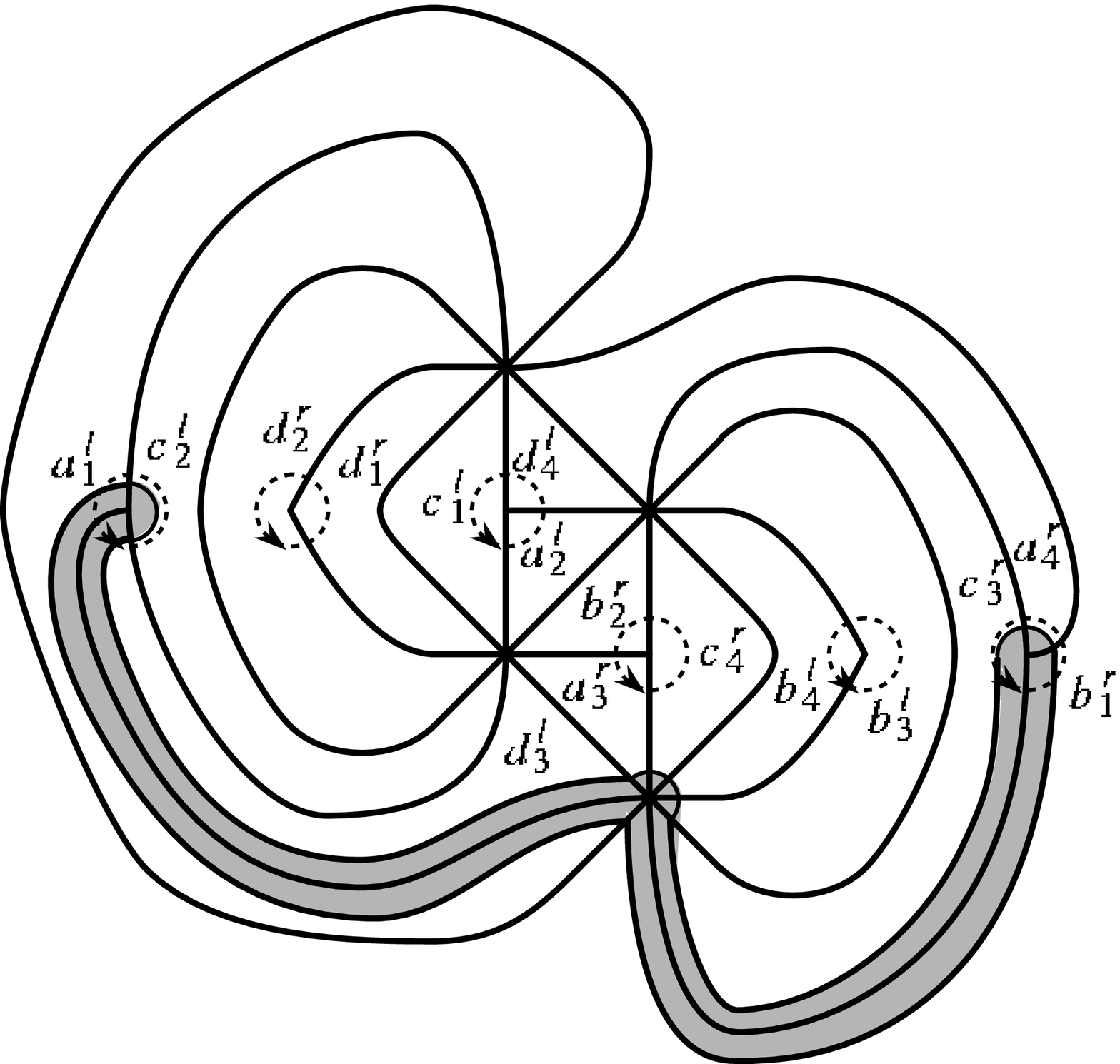}
\caption{The intersections of the regular neighborhood of Leaf~$\beta_1$
and the boundary of the bigger regular neighborhoods of ${\rm G}$ and
$\pm\infty$ are indicated by the shaded regions.}
\label{fig:shaded}
\end{figure}

Now we want to know the torus boundary of the regular neighborhood of this new
vertex.
To do that we first look at the regular neighborhood of Leaf~$\beta_1$
(= Leaf~$\delta_3$).
Figures~\ref{fig:nbd_octa1} and \ref{fig:nbd_octa4} indicate how it appears
in the first and the fourth octahedra respectively and how it is collapsed
if we collapse the leaf to a point.
In the tetrahedron decomposition of $S^3\setminus(K\cup\pm\infty)$ the
intersections of the regular neighborhood of Leaf~$\beta_1$ and the boundaries
of the (bigger) regular neighborhoods of ${\rm G}$ and $\pm\infty$ look like the
shaded regions of Figure~\ref{fig:shaded}.
\par
Then we remove the shaded regions and glue the resulting boundaries following
the recipe in Figures~\ref{fig:nbd_octa1} and
\ref{fig:nbd_octa4}, i.e.
glue
\begin{itemize}
  \item $c_1^+$ and $c_1^r$, $d_1^+$ and $d_1^l$, $a_3^-$ and $a_3^r$,
        $b_3^-$ and $b_3^l$, $a_4^r$ and $a_4^l$, $b_4^+$ and $b_4^l$,
        $c_4^+$ and $c_4^r$, $c_2^r$ and $c_2^l$, $d_2^-$ and $d_2^l$,
        and $a_2^-$ and $a_2^r$ with quadrangles, and
  \item $\{a_1^+, a_1^r, a_1^l\}$, $\{b_1^+, b_1^r, b_1^l\}$,
        $\{c_3^-, c_3^r, c_3^l\}$, $\{d_3^-, d_3^r, d_3^l\}$
        with hexagons.
\end{itemize}
\begin{figure}[h]
\includegraphics[scale=0.45]{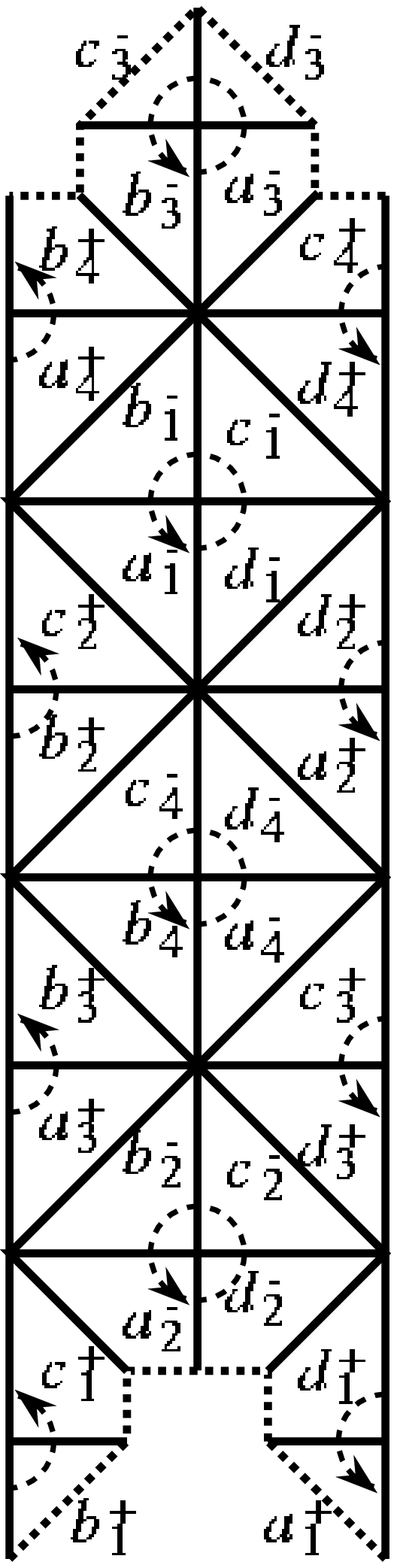}\\
\includegraphics[scale=0.45]{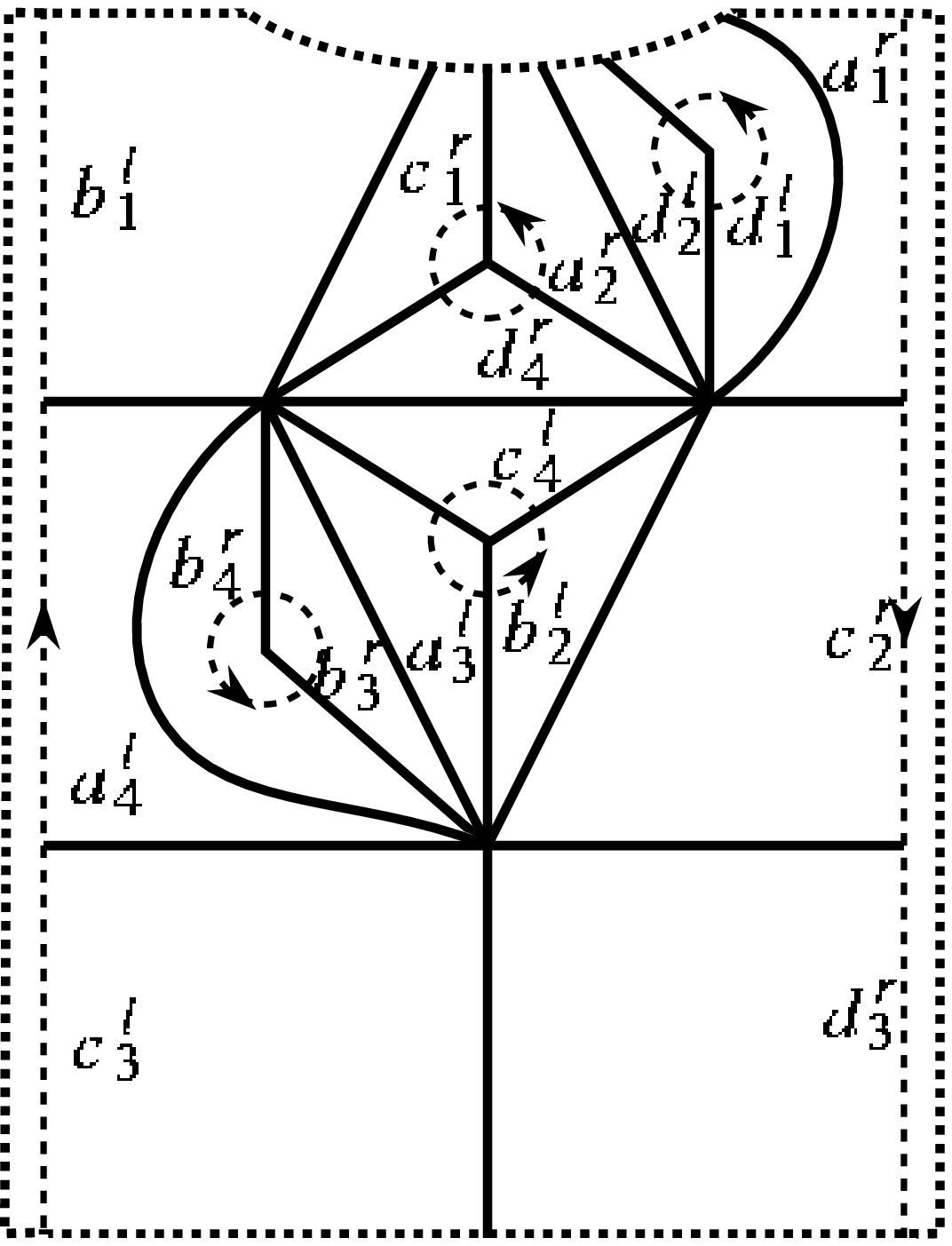}\qquad
\includegraphics[scale=0.45]{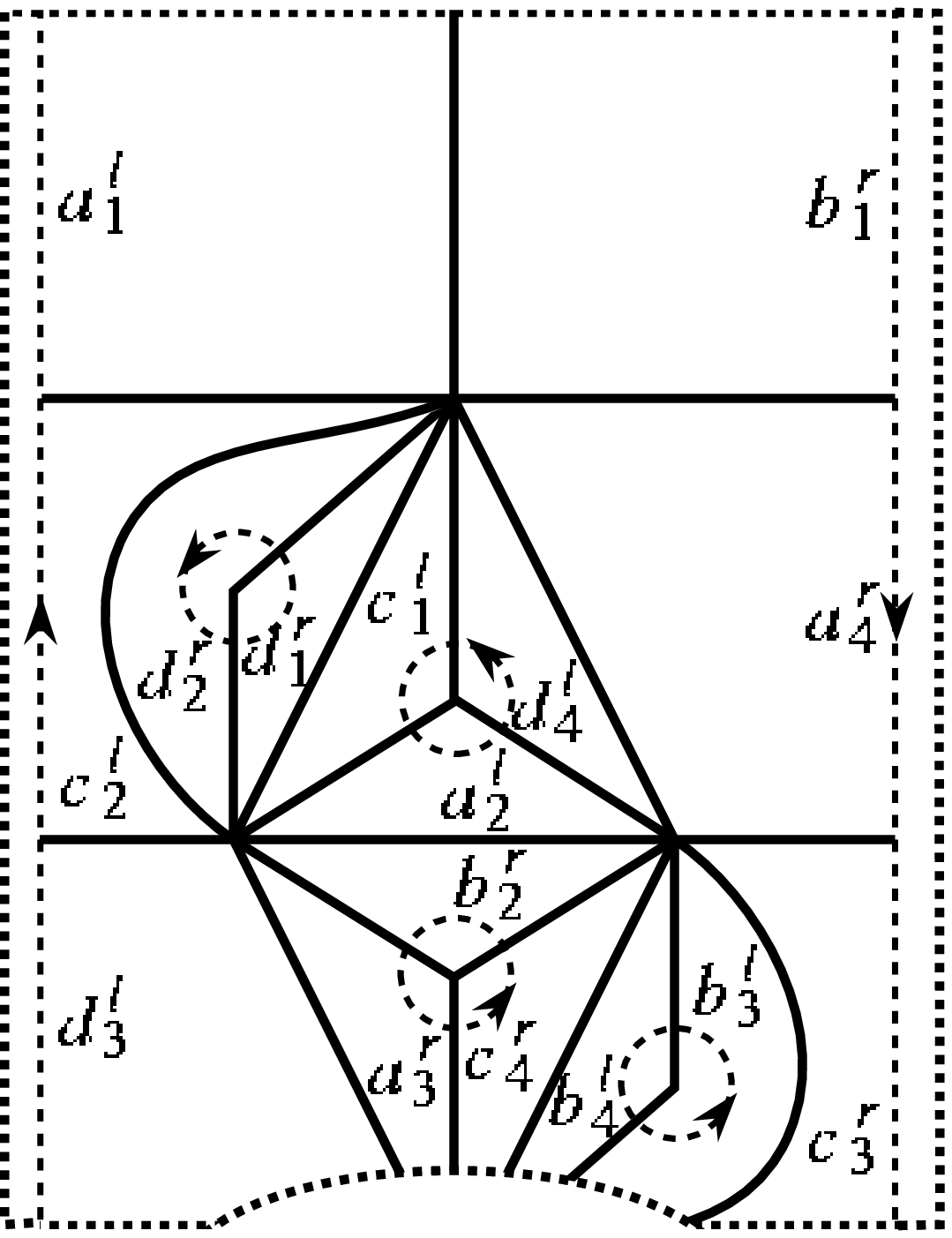}
\caption{Remove the shaded regions in Figure~\ref{fig:shaded} and deform the
resulting pictures.}
\label{fig:cut_infinity}
\end{figure}
Then we have the following pictures of the boundary of the regular neighborhood
of $\beta_1$ (Figures~\ref{fig:cut_infinity} and \ref{fig:cusp}).
\begin{figure}[h]
\includegraphics[scale=0.6]{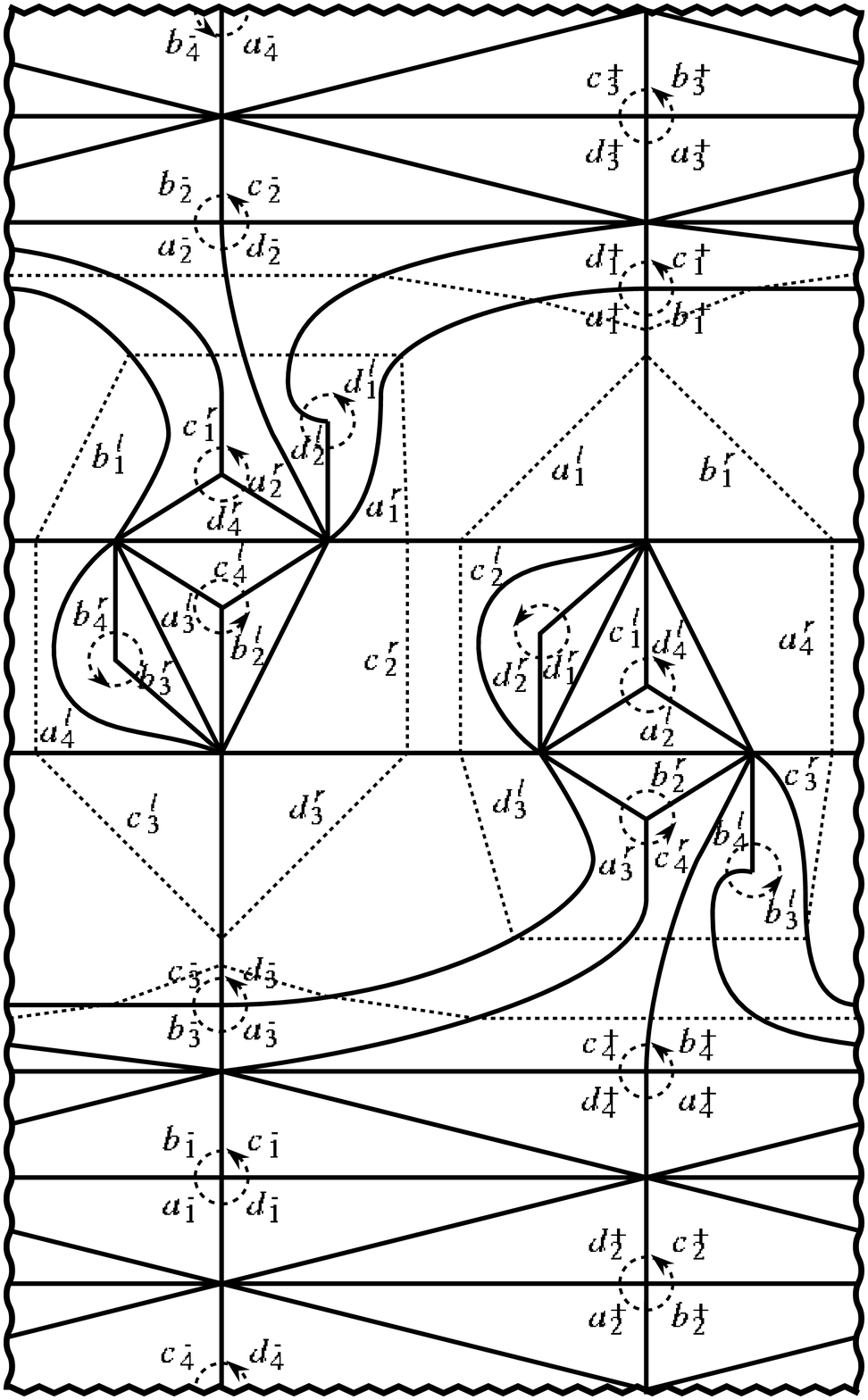}
\caption{The boundary of the regular neighborhood of $\beta_1$.}
\label{fig:cusp}
\end{figure}
To collapse it we again follow the recipe in Figures~\ref{fig:nbd_octa1}
and \ref{fig:nbd_octa4} and obtain the picture of the boundary
of the regular neighborhood of the collapsed $\beta_1$.
In Figure~$\ref{fig:collapsed_cusp}$ the edges with same numbers are identified.
\begin{figure}[h]
\includegraphics[scale=0.6]{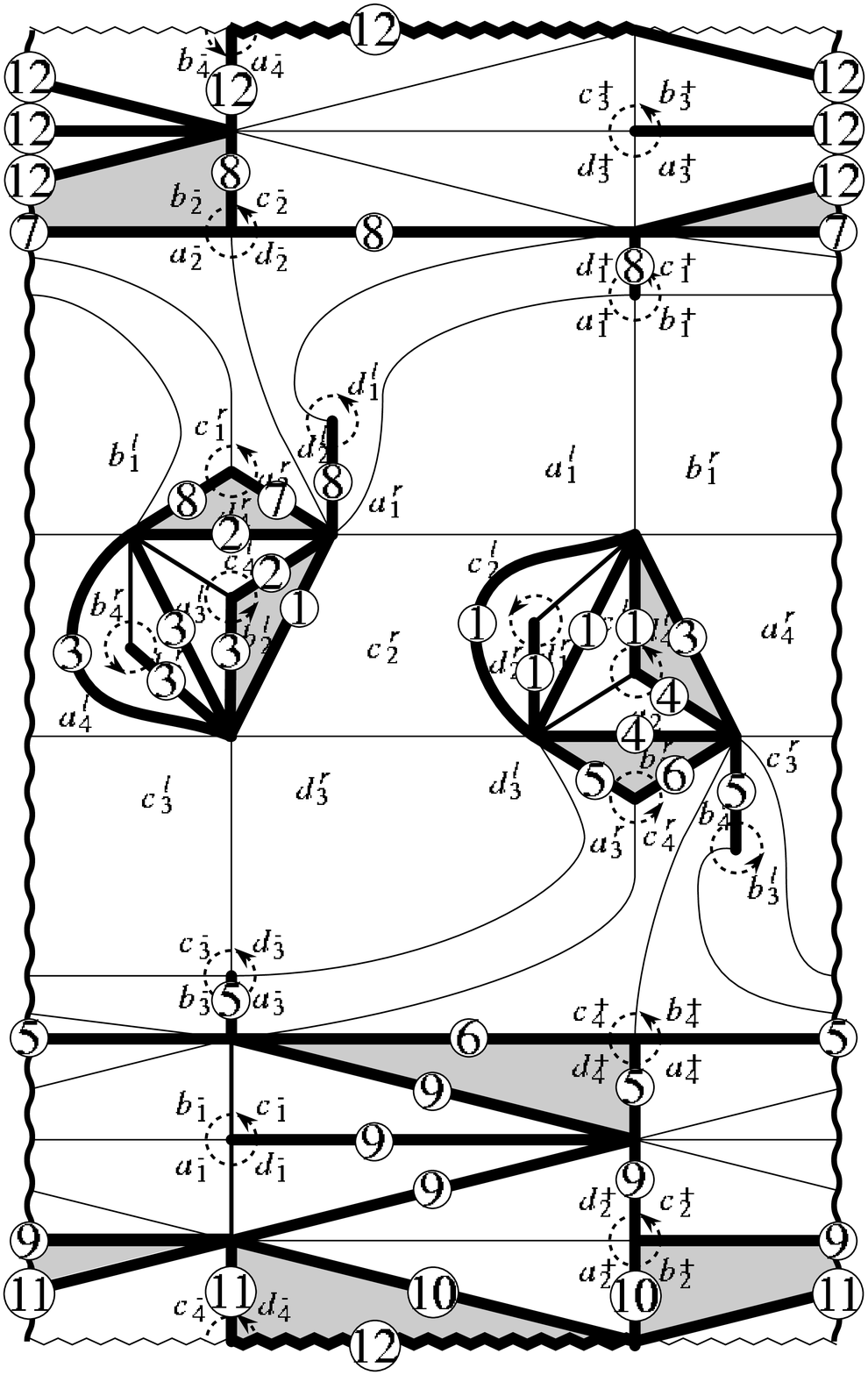}
\caption{The boundary of the regular neighborhood of $\beta_1$ after it
is collapsed to a point.
The edges with same numbers are identified and the shaded regions will
survive.}
\label{fig:collapsed_cusp}
\end{figure}
We finally have the triangulation of the cusp torus of $S^3\setminus{K}$
(Figure~\ref{fig:cusp_torus}).
In the figure $\lambda$ and $\mu$ denote the latitude (longitude) and the
meridian of $K$ respectively.
\begin{figure}[h]
\includegraphics[scale=0.6]{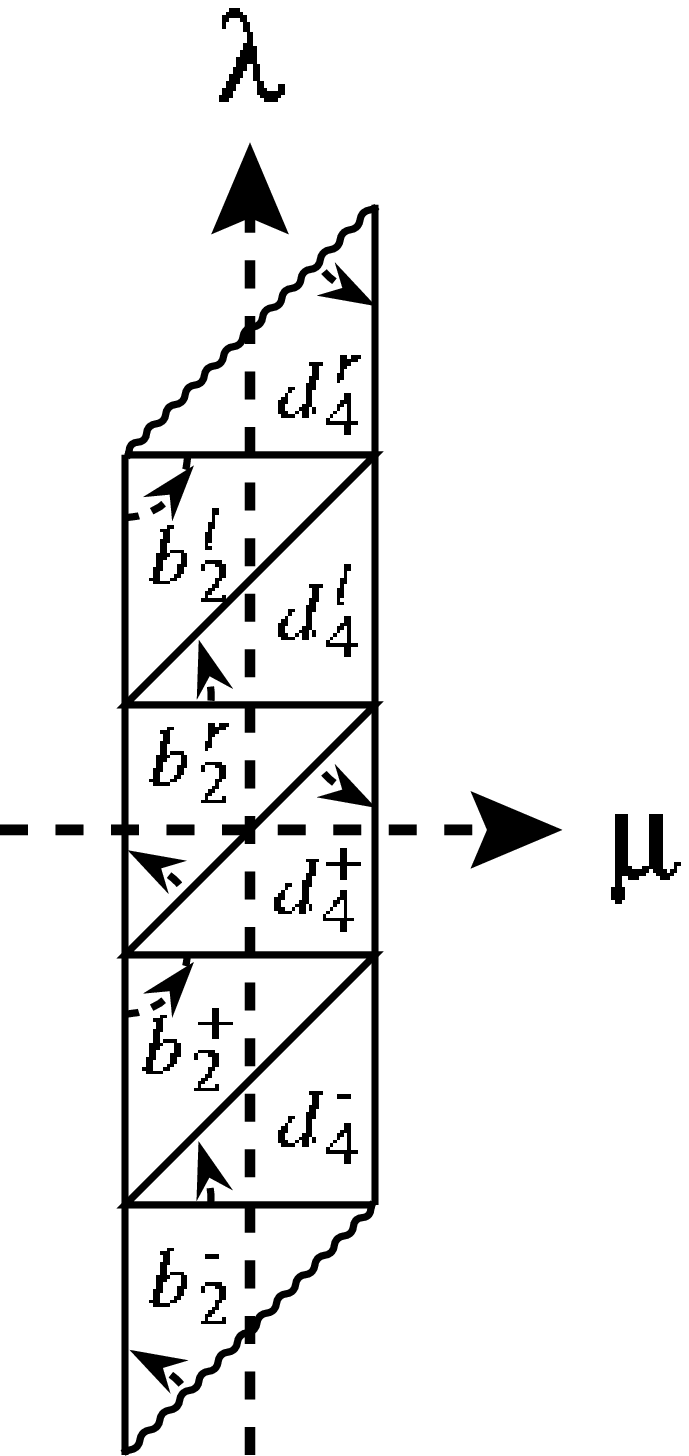}
\caption{The cusp torus of $S^3\setminus{K}$ with the latitude $\lambda$ and
the meridian $\mu$.}
\label{fig:cusp_torus}
\end{figure}
\subsection{Hyperbolicity condition}
\begin{figure}[h]
\includegraphics[scale=0.45]{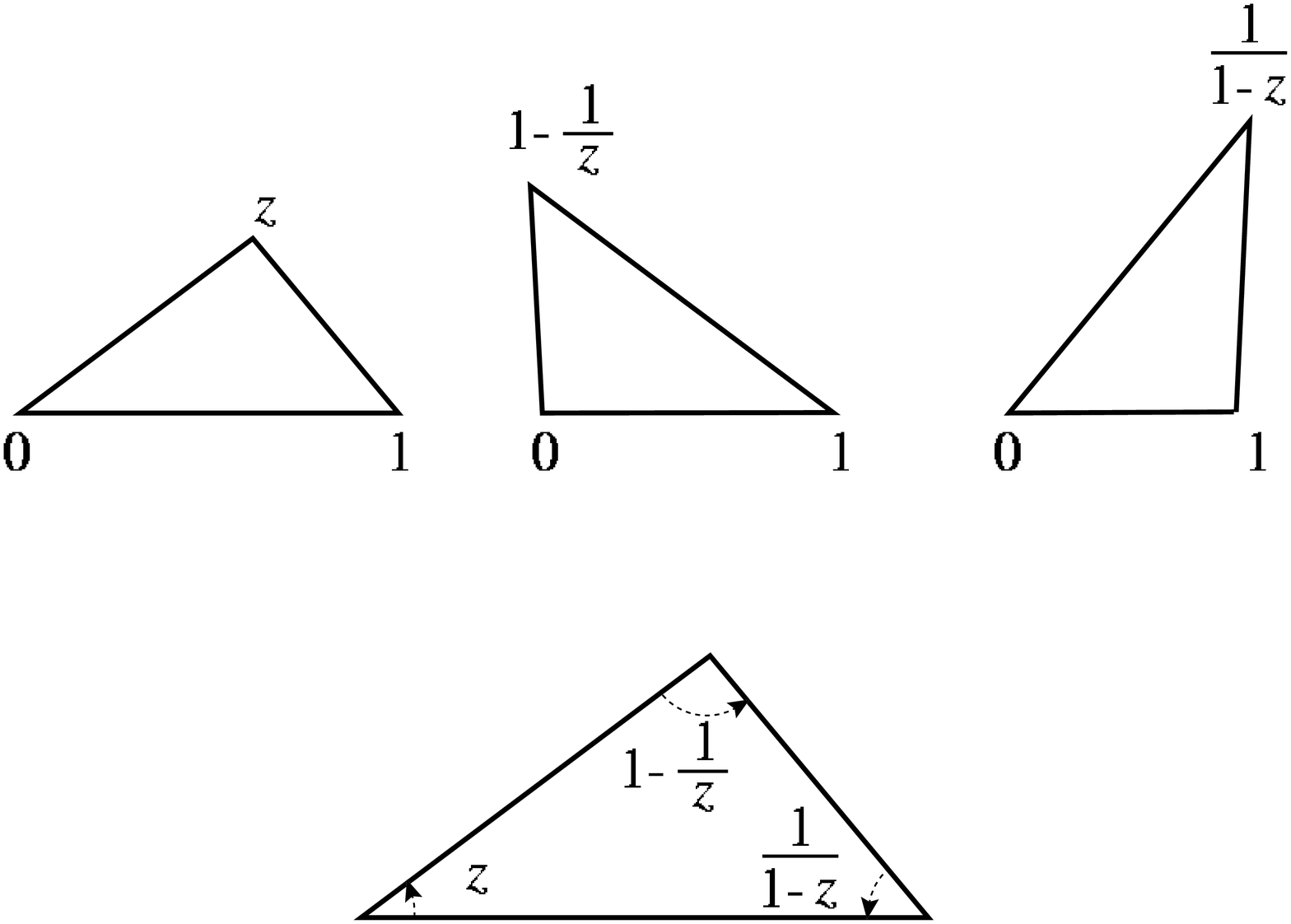}
\caption{Similar triangles parameterized by different complex numbers.}
\label{fig:triangle}
\end{figure}

Now we have a {\em topological} ideal tetrahedron decomposition of
$S^3\setminus{K}$.
To introduce a complete hyperbolic structure, we regard each tetrahedron
as an ideal hyperbolic tetrahedron parameterized by a complex number
with positive imaginary part.
This complex parameter $z$ is chosen as follows.
First recall that an ideal hyperbolic tetrahedron can be characterized by
a Euclidean triangle.
(If we consider the upper half space model of the hyperbolic space, then an
ideal hyperbolic tetrahedron can be located so that one ideal vertex is at the
infinity and the other three are on the $xy$-plane.
Then one face of the tetrahedron is in a hemisphere and the other three are in
upper half planes perpendicular to the $xy$-plane.
Our triangle is the intersection of the three half planes and the $xy$-plane.)
If one puts the triangle on the complex plane with two vertices on $0$ and $1$,
then the other is on $z$.
Note that there are three ways to choose such a parameter
(Figure~\ref{fig:triangle}).
For the tetrahedron with vertex $b_2^{\ast}$ we associate
$b\in\C$ and that with $d_4^{\ast}$ we associate
$d\in\C$, where $\ast$ is $+$, $-$, $r$, or $l$.
In the cusp torus the triangle labeled $b_2^{\ast}$ and $d_4^{\ast}$
are similar to the triangles parameterized by $b$ and $d$ respectively
for any $\ast$.
\par
It is well known that if the following two conditions are satisfied,
then an ideal tetrahedron decomposition parameterized as above defines
a hyperbolic structure (more precisely complete hyperbolic structure).
(See for example \cite{Neumann/Zagier:TOPOL85,Neumann92}.)
\begin{enumerate}
  \item Consistency condition: at each vertex in the cusp torus, the product of
  the parameter around the vertex is $1$, and
  \item Cusp condition: the products of the parameters along the latitude
  and the meridian are $1$.
  Here the parameter of a triangle along a loop is the parameter associated to
  the small triangle cut by the loop.
\end{enumerate}
\begin{rem}
The consistency condition means that at each edge of the tetrahedron
decomposition the sum of the dihedral angles is equal to $2\pi$.
The cusp condition means that the cusp torus is Euclidean.
\end{rem}
\par
In our example we have two edges, one connecting $b_2^l=b_2^-$ and
$b_2^r=b_2^+$
and the other connecting $d_4^l=d_4^+$ and $d_4^r=d_4^-$.
So from the consistency condition we have the following two equations
since the consistency conditions of two vertices on the torus are the
same if they belong to the same edge.
\begin{align*}
  b\left(\frac{1}{1-d}\right)b\left(1-\frac{1}{d}\right)
  \left(1-\frac{1}{b}\right)\left(1-\frac{1}{d}\right)&=1,
  \\
  d\left(\frac{1}{1-b}\right)\left(\frac{1}{1-d}\right)
  \left(\frac{1}{1-b}\right)d\left(1-\frac{1}{b}\right)&=1.
\end{align*}
From the cusp condition we have the following two equations.
\begin{align}
  \mu&: b^{-1}d=1,
  \label{eq:meridian}
  \\
  \lambda&:
  \left(\frac{1}{1-b}\right)^{-1}d
  \left(1-\frac{1}{b}\right)^{-1}\left(\frac{1}{1-d}\right)
  \\
  &\quad\times\left(\frac{1}{1-b}\right)^{-1}d
  \left(1-\frac{1}{b}\right)^{-1}\left(\frac{1}{1-d}\right)
  =1.\notag
\end{align}
From \eqref{eq:meridian} we have $b=d$ and then the rest are all equivalent
to $b^2-b+1=0$.
Since the imaginary part of $b$ should be positive, we have
$b=d=\exp(\pi\sqrt{-1}/3)$ and so the two tetrahedra should be regular.
\section{Analysis, and Harmony}\label{sec:harmony}
In \S\ref{sec:algebra} we know that if we calculate the Kashaev invariant
of the figure-eight knot only the parameter $j$ contributes.
Due to Kashaev \cite{Kashaev:LETMP97}, the asymptotic behavior of
$\langle 4_1\rangle_N$ for large $N$ is
\begin{equation*}
  \int\exp
  \left[
    \frac{N}{2\pi\sqrt{-1}}
    \left\{
      -\operatorname{Li}_2(z)+\operatorname{Li}_2(z^{-1})
    \right\}
  \right]\,dz,
\end{equation*}
where $\operatorname{Li}_2(z)=-\int_{0}^{z}\log(1-u)/u\,du$ is
the dilogarithm function.
Note that $z$ corresponds to $q^{j}$ and that $-\operatorname{Li}_2(z)$
and $\operatorname{Li}_2(z^{-1})$
corresponds to $(q)_{j}$ and $(q^{-1})_{j}$ respectively.
We can apply the saddle point method to this integral and know that
for large $N$ it is approximated by the value
\begin{equation*}
  \exp\left[\frac{N\,V(z_0)}{2\pi\sqrt{-1}}\right]
\end{equation*}
with $V(z):=-\operatorname{Li}_2(z)+\operatorname{Li}_2(z^{-1})$
and $z_0$ satisfies $d\,V(z_0)/d\,z=0$.
Since
\begin{equation*}
  \frac{d\,V(z)}{d\,z}=\frac{\log(1-z)}{z}+\frac{\log(1-z^{-1})}{z},
\end{equation*}
$z_0$ satisfies
\begin{equation}\label{eq:saddle}
  z_0^2-z_0+1=0.
\end{equation}
Now the Kashaev conjecture follows in this case since $V(z_0)$ is the
volume of the figure-eight knot complement, which is twice the volume 
of the regular tetrahedron.
\par
On the other hand in \S\ref{sec:geometry} we know that only two tetrahedra
parameterized by $b$ and $d$ survive.
Moreover from the cusp condition these two parameters coincide and satisfy
the equation
\begin{equation*}
  b^2-b+1=0,
\end{equation*}
which coincides with \eqref{eq:saddle}!
Observe that $-\operatorname{Li}_2(z)$ corresponds to $(q)_{j-1}=N/(q)_{[-j]}$
which comes from the angle between the edges around the second crossing labeled
$0$ and $j$, and $\operatorname{Li}_2(z^{-1})$ to
$(q^{-1})_{j-1}=N/(q^{-1})_{[-j]}$ which comes from the angle around the fourth
crossing between the edges labeled $0$ and $j$.
\par
Yokota proved in \cite{Yokota:volume} similar situations occur
for much more general knots including non-alternating knots, showing
that for a fixed knot diagram the angles between the surviving labelings of the
summation in $\langle{K}\rangle_N$ correspond to the surviving dihedral angles
in the ideal tetrahedron decomposition and that the partial differential
equations defining the saddle point of the integral approximating
$\langle{K}\rangle_N$ for large $N$ coincide with the hyperbolicity condition.
I hope that this note will help the reader to read Yokota's rather complicated
but beautiful paper \cite{Yokota:volume}.
\bibliography{mrabbrev,hitoshi}

\ifx\undefined\bysame
\newcommand{\bysame}{\leavevmode\hbox to3em{\hrulefill}\,}
\fi
\begin{thebibliography}{10}

\bibitem{Kashaev:MODPLA95}
R.M. Kashaev, {\em A link invariant from quantum dilogarithm}, Modern Phys.
  Lett. A {\bf 10} (1995), no.~19, 1409--1418.

\bibitem{Kashaev:LETMP97}
\bysame, {\em The hyperbolic volume of knots from the quantum dilogarithm},
  Lett. Math. Phys. {\bf 39} (1997), no.~3, 269--275.

\bibitem{Murakami/Murakami:volume}
H.~Murakami and J.~Murakami, {\em The colored {J}ones polynomials and the
  simplicial volume of a knot}, Report No. 20, 1998/99, Institut
  Mittag-Leffler, arXiv:math.GT/9905075.

\bibitem{Murakami/Murakami/Okamoto/Takata/Yokota:CS}
H.~Murakami, J.~Murakami, M.~Okamoto, T.~Takata, and Y.~Yokota, {\em Kashaev's
  conjecture and the {C}hern--{S}imons invariants of knots and links}, in
  preparation.

\bibitem{Neumann92}
W.D. Neumann, {\em Combinatorics of triangulations and the {C}hern-{S}imons
  invariant for hyperbolic $3$-manifolds}, Topology '90 (Columbus, OH, 1990)
  (Berlin), de Gruyter, Berlin, 1992, pp.~243--271.

\bibitem{Neumann/Zagier:TOPOL85}
W.D. Neumann and D.~Zagier, {\em Volumes of hyperbolic three-manifolds},
  Topology {\bf 24} (1985), no.~3, 307--332.

\bibitem{D.Thurston:Grenoble}
D.~Thurston, {\em Hyperbolic volume and the {J}ones polynomial}, Lecture notes,
  {\'E}cole d'{\'e}t{\'e} de Math{\'e}matiques `Invariants de n{\oe}uds et de
  vari{\'e}t{\'e}s de dimension $3$', Institut Fourier - UMR 5582 du CNRS et de
  l'UJF Grenoble (France) du 21 juin au 9 juillet 1999.

\bibitem{Turaev:INVEM88}
V.~G. Turaev, {\em The {Y}ang-{B}axter equation and invariants of links},
  Invent. Math. {\bf 92} (1988), 527--553.

\bibitem{Yokota:volume}
Y.~Yokota, {\em On the volume conjecture for hyperbolic knots}, in preparation.

\bibitem{Yokota:Topology_Symposium2000}
\bysame, {\em On the volume conjecture for hyperbolic knots}, Proceedings of
  the 47th Topology Symposium (Inamori Hall, Kagoshima University), July 2000,
  pp.~38--44.

\bibitem{Yokota:Murasugi70}
\bysame, {\em On the volume conjecture of hyperbolic knots}, Knot Theory --
  dedicated to Professor Kunio Murasugi for his 70th birthday (M.~Sakuma, ed.),
  March 2000, pp.~362--367.

\end{thebibliography}
\bibliographystyle{amsplain}
\end{document}